# PHANTOM MAPS, SNT-THEORY, AND NATURAL FILTRATIONS ON $\varprojlim^1$ SETS

PIERRE GHIENNE

Abstract. We study the so-called Gray filtration on the set of phantom maps between two spaces. Using both its algebraic characterization and the Sullivan completion approach to phantom maps, we generalize some of the recent results of Le, McGibbon and Strom. We particularly emphasize on the set of phantom maps with infinite Gray index, describing it in an original algebraic way.

We furthermore introduce and study a natural filtration on SNT-sets (that is sets of homotopy types of spaces having the same $n$-type for all $n$), which appears to have the same algebraic characterization of the Gray one on phantom maps. For spaces whose rational homotopy type is that of an $H$-space or a co-$H$-space, we establish criteria permitting to determinate those subsets of this filtration which are non trivial, generalizing work of McGibbon and Møller.

We finally describe algebraically the natural connection between phantom maps and SNT-theory, associating to a phantom map its homotopy fiber or cofiber. We use this description to show that this connection respect filtrations, and to find generic examples of spaces for which the filtration on the corresponding SNT-set consists of infinitely many strict inclusions.

## 1. Introduction

In this paper we propose to study natural filtrations on the set of phantom maps between two spaces, as well as on sets of homotopy types of spaces having the same $n$-type for all $n$.

Unless specifically precised, all spaces have the homotopy type of CW-complexes, are pointed, and maps between them are pointed maps.

Recall that a map $f : X \to Y$ is said to be a $\mathcal{F}$-phantom map, respectively a $\mathcal{F}D$-phantom map, if, for any CW-complex $K$ which is finite, respectively finite dimensional, and for any map $g : K \to X$, the composition $f \circ g : K \to Y$ is inessential.

Any $\mathcal{F}D$-phantom map is in particular a $\mathcal{F}$-phantom map, but the converse is false in general (see the survey papers of McGibbon [12] or Roitberg [20] for examples). We write $\mathrm{Ph}_{\mathcal{F}}(X,Y)$, respectively $\mathrm{Ph}_{\mathcal{F}D}(X,Y)$, for the subset of $[X,Y]$ consisting of homotopy classes of $\mathcal{F}$-phantom maps, respectively $\mathcal{F}D$-phantom maps. If $X$ has finite type, which means that $X$ admits finite skeleta, then the two notions of phantom maps do coincide, and we can write simply $\mathrm{Ph}(X,Y)$.

The filtration on $\mathcal{F}D$-phantom maps we work with, namely the *Gray filtration*, goes back to [4], and has been recently studied in [10] and [18]. In what follows we generalize the Gray filtration, and some of the results of [10], [18], to the case of $\mathcal{F}$-phantom maps. We also emphasize on and describe algebraically the set of $\mathcal{F}$-phantom maps with infinite *Gray index*.

Recall now that for a given space $X$, we denote $\mathrm{SNT}(X)$ the set of homotopy types of spaces $Y$ having the same $n$-type of $X$ for all $n$. That is to say that $X^{(n)}$ and $Y^{(n)}$, the Postnikov sections of $X$ and $Y$ through dimension $n$, have the same homotopy type for all $n$.

We introduce and study a natural filtration on this set, which appears to have the same algebraic characterization of the Gray one on $\mathcal{F}D$-phantom maps. Working algebraically, we show that the natural connection between $\mathcal{F}D$-phantom maps and SNT-theory (to a phantom map we associate its homotopy fiber or its homotopy cofiber) respect these filtrations. We use this result to find generic examples of spaces $X$ for which the filtration on $\mathrm{SNT}(X)$ consists of infinitely many strict inclusions.

We state our main results in Sections 2 and 3. The rest of the paper is devoted to proofs and complementary results.

## 2. The Gray filtration on phantom maps

The definition of the *Gray index* of a phantom map we give below deserves some comments. Originally [4], the Gray index is only define for $\mathcal{F}D$-phantom maps, by means of a particular cell-decomposition of the domain $X$. Our Definition 2.1 coincide with the usual one for $\mathcal{F}D$-phantom maps [10], and has the advantage to generalize naturally to the case of $\mathcal{F}$-phantom maps.

**Definition 2.1.** *Let $f : X \to Y$ be a $\mathcal{F}$-phantom map, respectively a $\mathcal{F}D$-phantom map. We say that the* Gray index *of $f$ is greater than $k$, and write $G(f) \geq k$, if there exists a lift $\tilde{f}$ of $f$ through $Y\langle k\rangle$, the $k$-connected cover of $Y$,*

$$\begin{array}{ccc} & & Y\langle k\rangle \\ & \overset{\tilde{f}}{\nearrow} & \downarrow \\ X & \xrightarrow[f]{} & Y \end{array}$$

*such that $\tilde{f}$ itself is a $\mathcal{F}$-phantom map, respectively a $\mathcal{F}D$-phantom map.*

Remark however the following fact: the existence of a lift $\tilde{f}$ is guaranteed if $f$ is $\mathcal{F}D$-phantom, whereas it is not if $f$ is only supposed to be $\mathcal{F}$-phantom.

As an easy example, note that the Gray index of the trivial map is infinite: $G(*) = \infty$.

The *Gray filtration* on phantom maps is defined by the following rule: a phantom map $f$ in filtration $k$ has Gray index $G(f) \geq k$. More precisely, we define

$$\mathrm{Ph}^k_{\mathcal{F}D}(X,Y) := \{f \in \mathrm{Ph}_{\mathcal{F}D}(X,Y) \mid G(f) \geq k\}$$

$$\mathrm{Ph}^k_{\mathcal{F}}(X,Y) := \{f \in \mathrm{Ph}_{\mathcal{F}}(X,Y) \mid G(f) \geq k\}$$

and we have two filtrations:

$$\begin{array}{ccccccc} \mathrm{Ph}_{\mathcal{F}D}(X,Y) & \supseteq \cdots \supseteq & \mathrm{Ph}^k_{\mathcal{F}D}(X,Y) & \supseteq \cdots \supseteq & \mathrm{Ph}^\infty_{\mathcal{F}D}(X,Y) := \bigcap_k \mathrm{Ph}^k_{\mathcal{F}D}(X,Y) \\ \cap\vert & & \cap\vert & & \cap\vert \\ \mathrm{Ph}_{\mathcal{F}}(X,Y) & \supseteq \cdots \supseteq & \mathrm{Ph}^k_{\mathcal{F}}(X,Y) & \supseteq \cdots \supseteq & \mathrm{Ph}^\infty_{\mathcal{F}}(X,Y) := \bigcap_k \mathrm{Ph}^k_{\mathcal{F}}(X,Y) \end{array}$$

Notice that we have also an inclusion:

$$\mathrm{Ph}^\infty_{\mathcal{F}}(X,Y) \subseteq \mathrm{Ph}_{\mathcal{F}D}(X,Y)$$

Indeed, a map $f \in \mathrm{Ph}^{\infty}_{\mathcal{F}}(X,Y)$ admits in particular, for any $k$, a lift through $Y\langle k\rangle$, and the composition $X \xrightarrow{f} Y \to Y^{(k)}$ is then inessential. The inclusion follows by classical theory.

The following example, due to McGibbon-Strom, shows that the Gray filtration can be highly non trivial:

**Example 2.2.** [18] Let $X := \mathbb{C}P(\infty)$ and $Y := S^2 \vee S^2$. Then $\mathrm{Ph}(X,Y) = \mathrm{Ph}^2(X,Y)$, and for any $k \geq 1$, we have $\mathrm{Ph}^{2k-1}(X,Y) = \mathrm{Ph}^{2k}(X,Y)$ whereas the inclusion $\mathrm{Ph}^{2k}(X,Y) \supset \mathrm{Ph}^{2k+1}(X,Y)$ is strict. Moreover, $\mathrm{Ph}^{\infty}(X,Y) = \{*\}$.

We list in the following theorem the results of [10] which are relevant to our present work.

**Theorem 2.3.** [10] *Let $X$ and $Y$ be nilpotent spaces of finite type.*
*1) If $\mathrm{Ph}(X,Y) = \mathrm{Ph}^{\infty}(X,Y)$, then in fact $\mathrm{Ph}(X,Y) = \{*\}$.*
*2) If $H^k(X;\pi_{k+1}(Y)\otimes\mathbb{Q}) \cong 0$, then $\mathrm{Ph}^k(X,Y) = \mathrm{Ph}^{k+1}(X,Y)$.*
*3) If $H^m(X;\pi_{m+1}(Y)\otimes\mathbb{Q}) \cong 0$ for any $m \geq k$, then $\mathrm{Ph}^k(X,Y) = \{*\}$.*

The first result we wish to present is the following generalization of Theorem 2.3.

**Theorem 2.4.**
*1) Let $X$ and $Y$ be nilpotent spaces of finite type. If $\mathrm{Ph}^k(X,Y) = \mathrm{Ph}^{\infty}(X,Y)$ for some $k$, then in fact $\mathrm{Ph}^k(X,Y) = \{*\}$.*
*2) Let $X$ be any 0-connected space, and $Y$ nilpotent of finite type.*
*If $H^k(X;\pi_{k+1}(Y)\otimes\mathbb{Q}) \cong 0$, then $\mathrm{Ph}^k_{\mathcal{F}}(X,Y) = \mathrm{Ph}^{k+1}_{\mathcal{F}}(X,Y)$.*
*3) Let $X$ be any 0-connected space, and $Y$ nilpotent of finite type.*
*If $H^m(X;\pi_{m+1}(Y)\otimes\mathbb{Q}) \cong 0$ for any $m \geq k$, then $\mathrm{Ph}^k_{\mathcal{F}}(X,Y) = \{*\}$.*

**Remark 2.5.** Point *1)* of Theorem 2.4 is more generally true for $\mathcal{F}D$-phantom maps between any spaces $X$ and $Y$ such that the groups $[X,\Omega Y^{(n)}]$ are countable for $n$ large enough.

Points *2)* and *3)* of Theorem 2.4 appears as natural generalization of the same Points in Theorem 2.3. Point *1)* of Theorem 2.4 generalizes both Points *1)* and *3)* of Theorem 2.3, saying that if the filtration stabilizes at any stage and for any reason, then it stabilizes to the trivial set.

We emphasize now on the sets $\mathrm{Ph}^{\infty}_{\mathcal{F}D}(X,Y)$ and $\mathrm{Ph}^{\infty}_{\mathcal{F}}(X,Y)$ of phantom maps with infinite Gray index. In his thesis [4], Gray claims that the only $\mathcal{F}D$-phantom map with infinite Gray index, between any spaces, is the trivial map. Counter-examples have been recently found by McGibbon-Strom [18]: to any space $X$, they associate canonically a $\mathcal{F}D$-phantom map $\Phi_X : X \to W(X)$ with infinite Gray index, and show this map is essential under some additional hypothesis on $X$.

However, for such $X$, even though these $X$ are finite type, the target $W(X)$ is definitively not. We quote the following conjecture arising from [10] and [18].

**Conjecture.** If $X$ and $Y$ are nilpotent of finite type, then $\mathrm{Ph}^{\infty}(X,Y) = \{*\}$.

Special cases of this conjecture appear already in the previous two theorems, as well as in the following result. Here $\widehat{Y}$ denotes the profinite completion of $Y$.

**Theorem 2.6.** [18] *Let $X$ and $Y$ be nilpotent of finite type, and suppose that the function space $\mathrm{Map}_*(X,\widehat{Y})$ is weakly contractible. Then $\mathrm{Ph}^{\infty}(X,Y) = \{*\}$.*

In the context of $\mathcal{F}$-phantom maps, we provide the following generalization of Theorem 2.6.

**Theorem 2.7.** *Let $X$ be any 0-connected space, and $Y$ nilpotent of finite type.*
*1) If $[X, \Omega\widehat{Y}]$ is finite, then* $\mathrm{Ph}^{\infty}_{\mathcal{F}D}(X,Y) = \{*\} = \mathrm{Ph}^{\infty}_{\mathcal{F}}(X,Y)$.
*2) If $[X, \Omega\widehat{Y}]$ is countable, then the equality* $\mathrm{Ph}^{k}_{\mathcal{F}}(X,Y) = \mathrm{Ph}^{k+1}_{\mathcal{F}}(X,Y)$ *holds if and only if* $H^k(X; \pi_{k+1}(Y) \otimes \mathbb{Q}) \cong 0$.

For the sake of completeness, we list in Proposition 5.5 well known cases where Theorem 2.7 do apply.

Point *2)* of Theorem 2.7 may be its most striking part: the determination of the Gray filtration then amounts to an easy rational computation. It gives for example an easy proof of Example 2.2.

Most of our results so far will be proved using the following characterization of $\mathcal{F}$-phantom maps, following directly from Sullivan [22, Theorem 3.2]. Let $\widehat{e} : Y \to \widehat{Y}$ be the profinite completion of the space $Y$. If $Y$ is nilpotent of finite type, then, for any space $X$,

$$\mathrm{Ph}_{\mathcal{F}}(X,Y) = \mathrm{Ker}\left(\widehat{e}_* : [X,Y] \to [X,\widehat{Y}]\right)$$

We have to interpret the Gray index in that context. For any $n \geq 1$, we define a space $P^nY$ as a pullback, together with a map $p^n : Y \to P^nY$, as shown in the following diagram:

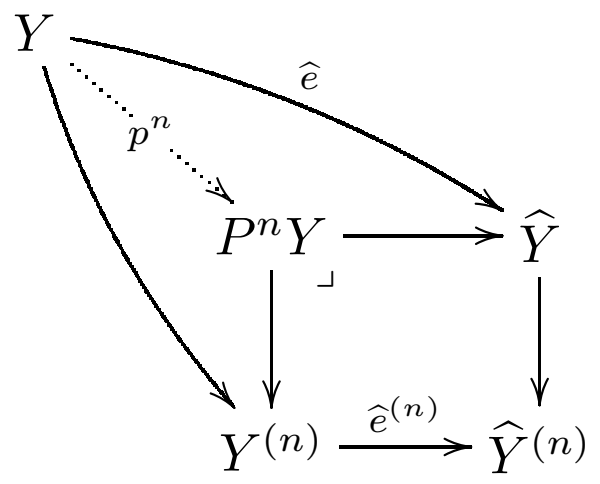

Notice that there might be several homotopy classes of maps $p^n : Y \to P^nY$ such that the above diagram commutes. Indeed, by classical theory, there is an action of $[Y, \Omega\widehat{Y}^{(n)}]$ on the set $[Y, P^nY]$, and if any map $p^n$ is chosen, all maps in the orbit of $p^n$ are also convenient. Working carefully, we can avoid this ambiguity (Corollary 5.10), and then construct an infinite tower approximating $Y$, that is to say that the maps $p^n$ induce an equivalence from $Y$ to the homotopy inverse limit of the tower.

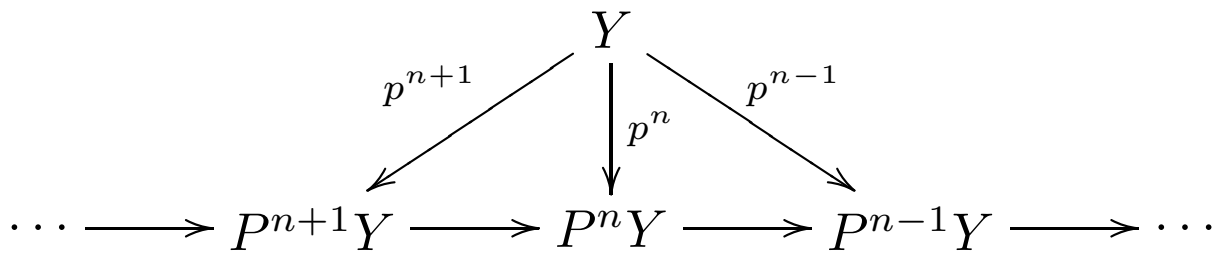

**Theorem 2.8.** *Let $Y$ be nilpotent of finite type, $X$ any space, and* $f \in \mathrm{Ph}_{\mathcal{F}}(X,Y)$.
*1) The Gray index $G(f) \geq n$ if and only if the composition* $X \xrightarrow{f} Y \xrightarrow{p^n} P^nY$ *is inessential.*
*2) We have a bijection* $\mathrm{Ph}^{\infty}_{\mathcal{F}}(X,Y) \cong \underset{\leftarrow n}{\lim}^1 [X, \Omega P^nY]$.

We refer to [1, p. 251] for the definition of the $\underset{\leftarrow}{\lim}^1$ of a tower of arbitrary groups. Recall from [1, p. 254] that we have a bijection:

$$\mathrm{Ph}_{\mathcal{F}D}(X,Y) \cong \underset{\leftarrow n}{\lim}^1 [X, \Omega Y^{(n)}]$$

Thanks to Point *2)* of Theorem 2.8, we can give an algebraic description of the inclusion $\mathrm{Ph}^{\infty}_{\mathcal{F}}(X,Y) \subseteq \mathrm{Ph}_{\mathcal{F}D}(X,Y)$: it is the $\underset{\leftarrow n}{\lim}^1$ of the homomorphisms $[X, \Omega P^nY] \to [X, \Omega Y^{(n)}]$ defined naturally by the maps $P^nY \to Y^{(n)}$.

For other $\varprojlim^1$ expressions of the set $\mathrm{Ph}^{\infty}_{\mathcal{F}}(X, Y)$, see Corollary 5.11.

The spaces $P^nY$ permit the following characterization of nilpotent spaces of finite type which can't be the target of a non-trivial $\mathcal{F}$-phantom map.

**Theorem 2.9.** *Let $Y$ be nilpotent of finite type, and $\overline{\Pi_{n\geq 1}\,\Omega P^nY}$ be a CW-approximation of the (possibly non CW-type) product $\Pi_{n\geq 1}\,\Omega P^nY$. Then $\mathrm{Ph}^{\infty}_{\mathcal{F}}(X, Y) = \{*\}$ for every space $X$ if and only if $\Omega Y$ is a retract of $\overline{\Pi_{n\geq 1}\,\Omega P^nY}$.*

It seems nevertheless difficult to see if a given space verify or not the hypothesis of Theorem 2.9. Here follows a first but incomplete result in this direction.

**Proposition 2.10.** *Let $Y$ be nilpotent of finite type. Suppose there exists a weak homotopy equivalence from $\Omega Y$ to some product $\Pi_\alpha L_\alpha$, in which each $L_\alpha$ is a finite type, rationally elliptic space. Then $\Omega Y$ is a retract of $\overline{\Pi_{n\geq 1}\,\Omega P^nY}$.*

**Example 2.11.** If $Y = \bigvee_{i\in I} S^{n_i}$ is any finite wedge of simply connected spheres, then for every space $X$ we have $\mathrm{Ph}^{\infty}_{\mathcal{F}}(X, Y) = \{*\}$. This follows from the above results and the Milnor-Moore Theorem.

We shall at present go back to the algebraic characterization of $\mathcal{F}D$-phantom maps. Frow now on and until we state Theorem 2.12 below, we follow Le-Strom [10].

Consider $\{G_n\}_n$ any tower of groups:

$$\cdots \to G_n \to G_{n-1} \to \cdots \to G_1 \to *$$

and define $G^n_k := \mathrm{Im}(G_n \to G_k)$ if $n \geq k$, and $G^n_k := G_n$ if $n < k$. We obtain a natural diagram of surjections between towers (left of the following picture), which induces a diagram of surjections between their $\varprojlim_n^1$ (right of the following picture):

$$\begin{array}{ccc} & & \{G^n_{k+1}\}_n \\ & \nearrow\!\!\!\nearrow & \Big\downarrow \\ \{G_n\}_n & & \\ & \searrow\!\!\!\searrow & \\ & & \{G^n_k\}_n \end{array} \qquad \Longrightarrow \qquad \begin{array}{ccc} & & \varprojlim_n^1 G^n_{k+1} \\ & \overset{p_{k+1}}{\nearrow\!\!\!\nearrow} & \Big\downarrow p_{k,k+1} \\ \varprojlim_n^1 G_n & & \\ & \underset{p_k}{\searrow\!\!\!\searrow} & \\ & & \varprojlim_n^1 G^n_k \end{array}$$

Set $L := \varprojlim_n^1 G_n$, and define $L^k := \mathrm{Ker}\, p_k$. We then have a filtration:

$$(\Delta) \qquad L = L^0 \supseteq L^1 \supseteq \cdots \supseteq L^k \supseteq \cdots \supseteq L^\infty := \bigcap_k L^k$$

We also define, for $m \geq k$, an equivalence relation $\sim_m$ on $L^k$. Two elements $x, y \in L^k$ are said to be $m$-equivalent, which we denote by $x \sim_m y$, if and only if they have the same image by $p_m : L^k \hookrightarrow L \to \varprojlim_n^1 G^n_m$.

**Theorem 2.12.** [10] *Let $G_n := [X, \Omega Y^{(n)}]$. The filtration $(\Delta)$ defined just above on $\mathrm{Ph}_{\mathcal{F}D}(X, Y) \cong L := \varprojlim_n^1 G_n$ is precisely the usual Gray filtration on $\mathcal{F}D$-phantom maps: for any $k$, we have $L^k \cong \mathrm{Ph}^k_{\mathcal{F}D}(X, Y)$.*

This result justify the terminology in the following definition.

**Definition 2.13.** *The* algebraic Gray filtration *on a $\varprojlim^1$ set $L := \varprojlim_n^1 G_n$ is the filtration $(\Delta)$.*

**Remark 2.14.** We draw attention to the following fact: it's possible to find towers $\{G_n\}_n$ and $\{G'_n\}_n$ admitting the same $\varprojlim^1$, but with radically different algebraic Gray filtration on it (see Proposition 5.13 or Example 6.2).

Notice that if $X$ and $Y$ are nilpotent spaces of finite type, then the groups $[X, \Omega Y^{(n)}]$ are nilpotent finitely generated [17]. Then Point *1)* of Theorem 2.4 follows from Theorem 2.12 and the following result.

**Theorem 2.15.** *Let $\{G_n\}_n$ be a tower of* countable *groups, and consider the algebraic Gray filtration on $L := \varprojlim_n^1 G_n$. If $L^\infty = L^k$ for some $k$, then in fact $L^k \cong \{*\}$.*

This result stated in all its generality will also have applications in SNT-theory.

## 3. A filtration on SNT($X$)

In this section we suppose spaces to be 0-connected. We denote by $\mathrm{Aut}(X)$ the group of homotopy classes of self-homotopy equivalences of a space $X$.

Recall from [24] that we have a bijection:

$$\mathrm{SNT}(X) \cong \varprojlim_n{}^1 \mathrm{Aut}(X^{(n)})$$

This fundamental result allowed a deep study of the subject by McGibbon-Møller [14, 15, 16], with in particular the calculation of $\mathrm{SNT}(X)$ for many spaces $X$. Nevertheless, little is known on the general structure of these sets, besides the following [14]: if $X$ is a nilpotent space with finite type over some subring of the rationals, then either $\mathrm{SNT}(X) = \{X\}$, or else it is uncountably large.

In this section we study the algebraic Gray filtration (Definition 2.13) inherited by $\mathrm{SNT}(X)$ from its $\varprojlim^1$ description:

$$\mathrm{SNT}(X) \supseteq \mathrm{SNT}^1(X) \supseteq \cdots \supseteq \mathrm{SNT}^k(X) \supseteq \cdots \supseteq \mathrm{SNT}^\infty(X) := \bigcap_k \mathrm{SNT}^k(X)$$

One may ask if this filtration could have been defined without any aid of the algebra, and indeed it could, as shown by the following result.

**Theorem 3.1.** *Let $X$ be a 0-connected space, and $k \geq 0$. Then $\mathrm{SNT}^k(X)$ is the set of homotopy types of spaces $Y$, such that there exists a collection of homotopy equivalences $f_n : X^{(n)} \to Y^{(n)}$, $n \geq 1$, which are compatible through range $k$. More precisely, for any couple $(n,k)$, the Postnikov sections $f_k^{(\min(n,k))}$ and $f_n^{(\min(n,k))}$ must be homotopy equivalent.*

Let us explain what is behind the proof of Theorem 3.1. Let $F$ be an homotopy functor from spaces to spaces. We the define a subset $\mathrm{SNT}_F(X)$ of $\mathrm{SNT}(X)$ by the following rule: the homotopy type of a space $Y$ belongs to $\mathrm{SNT}_F(X)$ if and only if there exists a collection of homotopy equivalences $f_n : X^{(n)} \to Y^{(n)}$, $n \geq 1$, which are compatible after applying the functor $F$. More precisely, for any $n \geq 1$, we must have a commutative diagram:

$$\begin{array}{ccc} F(X^{(n+1)}) & \xrightarrow{F(f_{n+1})} & F(Y^{(n+1)}) \\ \downarrow & & \downarrow \\ F(X^{(n)}) & \xrightarrow{F(f_n)} & F(Y^{(n)}) \end{array}$$

where vertical maps are obtained by applying $F$ to the canonical $X^{(n+1)} \to X^{(n)}$ and $Y^{(n+1)} \to Y^{(n)}$.

Let $\mathrm{Aut}_F(X)$ be the kernel of the group homomorphism $\mathrm{Aut}(X) \to \mathrm{Aut}(FX)$, and suppose now that $F$ commutes with Postnikov sections. Then the inclusions $\mathrm{Aut}_F(X^{(n)}) \hookrightarrow \mathrm{Aut}(X^{(n)})$ fit together to give a map:

$$j_F : \varprojlim_n{}^1 \mathrm{Aut}_F(X^{(n)}) \to \varprojlim_n{}^1 \mathrm{Aut}(X^{(n)}) \cong \mathrm{SNT}(X).$$

**Theorem 3.2.** *If $F$ commutes with Postnikov sections, then $\mathrm{SNT}_F(X)$ is exactly the image of the map $j_F$.*

We owe the idea of defining $\mathrm{SNT}_F(X)$ and the proof of Theorem 3.2 to a similar construction of Félix-Thomas [3]. The proof of Theorem 3.1 amounts to identify $\mathrm{SNT}^k(X)$ with $\mathrm{SNT}_F(X)$ when the functor $F$ is the Postnikov section $-^{(k)}$.

In complete analogy with the Gray filtration on phantom maps, our filtration on $\mathrm{SNT}(X)$ is endowed with the following stability properties. For a given space $X$, we denote $\mathrm{Aut}_k(X)$ the subgroup of $\mathrm{Aut}(X)$ of maps inducing the identity on $X^{(k)}$,

$$\mathrm{Aut}_k(X) := \mathrm{Ker}\left(\mathrm{Aut}(X) \to \mathrm{Aut}(X^{(k)})\right).$$

**Theorem 3.3.** *Let $X$ be a nilpotent space with finite type over some subring of the rationals.*
*1) If $\mathrm{SNT}^k(X) = \mathrm{SNT}^\infty(X)$ for some $k$, then in fact $\mathrm{SNT}^k(X) = \{X\}$.*
*2) If $\mathrm{Aut}_k(X^{(k+1)})$ is a finite group, then $\mathrm{SNT}^k(X) = \mathrm{SNT}^{k+1}(X)$. This holds in particular if $H^{k+1}(X; \pi_{k+1}(X) \otimes \mathbb{Q}) \cong 0$.*
*3) If $H^{m+1}(X; \pi_{m+1}(X) \otimes \mathbb{Q}) \cong 0$ for any $m \geq k$, then $\mathrm{SNT}^k(X) = \{X\}$.*

The direct computation of $\mathrm{SNT}(X)$ for a given space $X$ stays in general problematic, and it is certainly the same for its subsets $\mathrm{SNT}^k(X)$. We however offer the following criteria for $H_0$-spaces (whose rational homotopy type is an $H$-space) and co-$H_0$-spaces (whose rational homotopy type is a bouquet of spheres).

We clarify some notation. Let $\mathbb{Z}_\mathbb{P}$ be the integers localized at some set of primes $\mathbb{P}$, and let $X$ be a $\mathbb{P}$-local space. We denote by $\mathrm{Aut}_k(H^{\leq n}(X; \mathbb{Z}_\mathbb{P}))$ the subgroup of those ring automorphisms of the graded ring $H^{\leq n}(X; \mathbb{Z}_\mathbb{P})$, that preserve the degrees of homogeneous elements, and induce identity on $H^{\leq k}(X; \mathbb{Z}_\mathbb{P})$. Dually, $\mathrm{Aut}_k(\pi_{\leq n} X)$ denotes the group of those automorphisms of the graded $\mathbb{Z}_\mathbb{P}$-module $\pi_{\leq n} X$, that preserve the Whitehead product pairing, and induce identity on $\pi_{\leq k} X$.

**Theorem 3.4.**
*1) Let $X$ be a 1-connected $H_0$-space with finite type over $\mathbb{Z}_\mathbb{P}$. Then the following statements are equivalent:*
*(i) $\mathrm{SNT}^k(X) = \{X\}$;*
*(ii) For all $n \geq k$, the image of $\mathrm{Aut}_k(X) \to \mathrm{Aut}_k(X^{(n)})$ has finite index.*
*(iii) For all $n \geq k$, the image of $\mathrm{Aut}_k(X) \to \mathrm{Aut}_k(H^{\leq n}(X; \mathbb{Z}_\mathbb{P}))$ has finite index.*

*2) Let $X$ be a 1-connected co-$H_0$-space with finite type over $\mathbb{Z}_\mathbb{P}$. Then the following statements are equivalent:*
*(i) $\mathrm{SNT}^k(X) = \{X\}$;*
*(ii) For all $n \geq k$, the image of $\mathrm{Aut}_k(X) \to \mathrm{Aut}_k(X^{(n)})$ has finite index.*
*(iii) For all $n \geq k$, the image of $\mathrm{Aut}_k(X) \to \mathrm{Aut}_k(\pi_{\leq n} X)$ has finite index.*

Theorem 3.4 generalizes [14, Theorem 1] and [15, Theorem 1] which deal with the case $k = 0$. Our methods of proof are similar.

**Example 3.5.** The set $\mathrm{SNT}^{2m}(BSU(m))$ is trivial, whereas $\mathrm{SNT}^{2m-1}(BSU(m))$ is not. Similarly, $\mathrm{SNT}^{4m}(BSp(m))$ is trivial, whereas $\mathrm{SNT}^{4m-1}(BSp(m))$ is not.

We shall now use phantom map theory to find examples of spaces $Z$ for which the filtration on $\mathrm{SNT}(Z)$ consists of many strict inclusions. The connection between phantom maps and SNT-theory is not new (it goes back to [5]), but is particularly illuminating in the context of our work. It is given by the maps

$$\mathrm{SNT}(Y \vee \Sigma X) \xleftarrow{\mathrm{Cof}} \mathrm{Ph}_{\mathcal{F}D}(X, Y) \xrightarrow{\mathrm{Fib}} \mathrm{SNT}(X \times \Omega Y)$$

associating to a $\mathcal{F}D$-phantom map its homotopy cofiber or its homotopy fiber.

**Theorem 3.6.** *The maps* Fib *and* Cof *respect filtrations. More precisely, we have:*
*1) The image of* $\mathrm{Ph}^k_{\mathcal{F}D}(X,Y)$ *by* Fib *is included in* $\mathrm{SNT}^{k-1}(X \times \Omega Y)$.
*2) The image of* $\mathrm{Ph}^k_{\mathcal{F}D}(X,Y)$ *by* Cof *is included in* $\mathrm{SNT}^{k}(Y \vee \Sigma X)$.

The proof of Theorem 3.6 is based on an algebraic description of the maps Fib and Cof (Theorem 7.1).

Theorem 3.6 is the key point to prove the following examples. The Example 3.7 is a particular case of a more general result (Proposition 7.2). The Example 3.8 shows that, at least if we drop finiteness conditions on a space $Z$, then $\mathrm{SNT}^\infty(Z)$ need not be trivial.

**Example 3.7.** (Example 2.2 continued) There are infinitely many *strict* inclusions in the filtration on $\mathrm{SNT}\big(\mathbb{C}P(\infty) \times \Omega(S^2 \vee S^2)\big)$.

**Example 3.8.** Let $\Phi_X : X \to W(X)$ be the canonical $\mathcal{F}D$-phantom map with infinite Gray index out of $X$ [18]. Suppose $X$ is finite type, and that its cohomology is not locally finite as a module over the Steenrod Algebra for some prime $p$. Then the cofiber $W(X)/X$ of $\Phi_X$ belongs to $\mathrm{SNT}^\infty(W(X) \vee \Sigma X)$, and is not homotopy equivalent to $W(X) \vee \Sigma X$.

## 4. The algebraic approach to $\mathcal{F}D$-phantom maps

In this section we focus on the *algebraic Gray filtration* we defined on any $\varprojlim^1$ set (Definition 2.13). We prove Theorem 2.15 as well as some complementary results on $\mathcal{F}D$-phantom maps (Proposition 4.5).

Consider $\{G_n\}_n$ any tower of groups:

$$\cdots \to G_n \to G_{n-1} \to \cdots \to G_1 \to *$$

and define $G^n_k := \mathrm{Im}(G_n \to G_k)$ if $n \geq k$, and $G^n_k := G_n$ if $n < k$. We then obtain naturally, for each $k$, another tower $\{G^n_k\}_n$:

$$\cdots \hookrightarrow G^n_k \hookrightarrow G^{n-1}_k \hookrightarrow \cdots \hookrightarrow G^{k+1}_k \hookrightarrow G_k \to G_{k-1} \to \cdots \to G_1 \to *,$$

in which each map to the left of $G_k$ is an inclusion.

By definition, the tower $\{G_n\}_n$ is said to be *Mittag-Leffler* if and only if, for each $k$, the tower $\{G^n_k\}_n$ eventually stabilizes. Now, as almost all maps in the tower $\{G^n_k\}_n$ are inclusions, it is clear that $\{G^n_k\}_n$ stabilizes if and only if $\{G^n_k\}_n$ is itself Mittag-Leffler. This observation permits to reformulate [14, Theorem 2] as follows.

**Lemma 4.1.** *Let* $\{G_n\}_n$ *be a tower of countable groups. Then* $\varprojlim^1_n G_n \cong *$ *if and only if* $\varprojlim^1_n G^n_k \cong *$ *for any* $k$. *Otherwise,* $\varprojlim^1_n G_n$ *is an uncountable set.*

Remark now that for any $n$, and $m \geq k$, the map $G_m \to G_k$ induces a surjection $G^n_m \twoheadrightarrow G^n_k$. All these surjections, when $k$ and $m$ are fixed and $n$ varies, fit together to give a map of tower $\{G^n_m\}_n \to \{G^n_k\}_n$, inducing a surjection $\varprojlim^1_n G^n_m \twoheadrightarrow \varprojlim^1_n G^n_k$. This gives a right meaning to the following definition.

**Definition 4.2.**
*1) The* index $\mathrm{Ind}\{G_n\}_n$ *of the tower* $\{G_n\}_n$ *is the greatest* $k$ *such that* $\varprojlim^1_n G^n_k \cong *$. *If no such* $k$ *exists, we say that the index is infinite.*
*2) Fix* $k \geq 1$, *and set* $K^n_k := \mathrm{Ker}(G_n \twoheadrightarrow G^n_k)$. *The* k-th index $\mathrm{Ind}_k\{G_n\}_n$ *of the tower* $\{G_n\}_n$ *is the index of the tower* $\{K^n_k\}_n$.

**Remark 4.3.** If the groups $G_n$ are countable, Lemma 4.1 implies that $\varprojlim^1_n G_n \cong *$ if and only if $\mathrm{Ind}\{G_n\}_n = \infty$.

We shall now state and prove the key result of that section, which in particular describe the algebraic Gray filtration with respect to the indices defined above.

**Theorem 4.4.** *Let $\{G_n\}_n$ be a tower of countable groups, and consider $L^k$ the subsets of $L := \varprojlim_n{}^1 G_n$ defining the algebraic Gray filtration.*
*1) For any $k$, we have $L^k = *$ if and only if $\varprojlim_n{}^1 K_k^n = *$. Otherwise, $L_k$ is uncountable.*
*2) For any $m \geq k$, the set of $m$-equivalence classes in $L^k$, namely $L^k/_{\sim_m}$, is eiher $*$, or else uncountable. Moreover, if $\mathrm{Ker}(G_m \to G_k)$ is finite, then $L^k/_{\sim_m} = *$.*
*3) For any $m \geq k$, the equality $L^k = L^m$ holds if and only if $m \leq \mathrm{Ind}_k\{G_n\}_n$.*

*Proof.* For any $k$, and $n \geq m$, the map $G_n \to G_m$ induces another map $K_k^n \to K_k^m$. Set $K_{k,m}^n := \mathrm{Im}(K_k^n \to K_k^m)$ if $n \geq m$, and set $K_{k,m}^n := K_k^n$ if $n < m$. By definition, $\mathrm{Ind}_k\{G_n\}_n$ is the greatest $m$ such that $\varprojlim_n{}^1 K_{k,m}^n \cong *$.

Fix integers $m \geq k$. We claim that there is an action of the group $\varprojlim_n G_k^n$ on $\varprojlim_n{}^1 K_k^n$, as well as on $\varprojlim_n{}^1 K_{k,m}^n$, such that:

$$L^k \cong \varprojlim_n{}^1 K_k^n \,/\, \varprojlim_n G_k^n \qquad \text{and} \qquad L^k/_{\sim_m} \cong \varprojlim_n{}^1 K_{k,m}^n \,/\, \varprojlim_n G_k^n.$$

Let's assume that claim. From Lemma 4.1 we know that $\varprojlim_n{}^1 K_k^n$ and $\varprojlim_n{}^1 K_{k,m}^n$ are either $*$, or else uncountable. We then see the same must be true then for $L^k$ and $L^k/_{\sim_m}$. Indeed, if not $*$, these sets are quotient of an uncountable set by the countable group $\varprojlim_n G_k^n \subset G_k$. For the same reason, $L^k/_{\sim_m} = *$ if and only if $\varprojlim_n{}^1 K_{k,m}^n = *$, that is to say if and only if $m \leq \mathrm{Ind}_k\{G_n\}_n$. Clearly, $L^k/_{\sim_m} = *$ if and only if $L^k = L^m$. Moreover, for $n \geq m$, we have $K_{k,m}^n \subseteq K_k^m = \mathrm{Ker}(G_m \to G_k)$, and if this kernel is finite, we have clearly $\varprojlim_n{}^1 K_{k,m}^n = *$. The theorem follows.

We now show our claim. Consider the following *first* diagram:

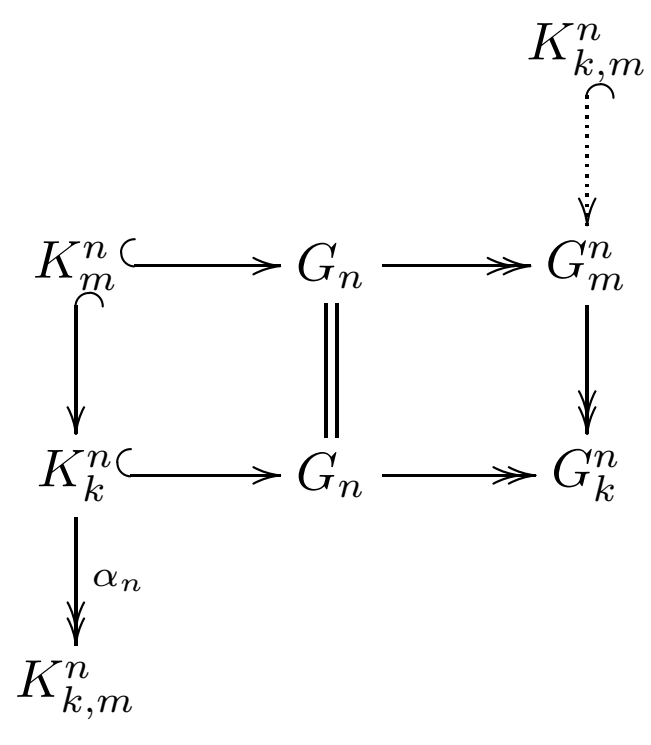

where the rows are exact by definition, and $\alpha_n$ is the obvious map. The left column is easily seen to be exact, and the snake lemma induces the dotted map such that the right column is also exact.

Consider now the following *second* commutative diagram:

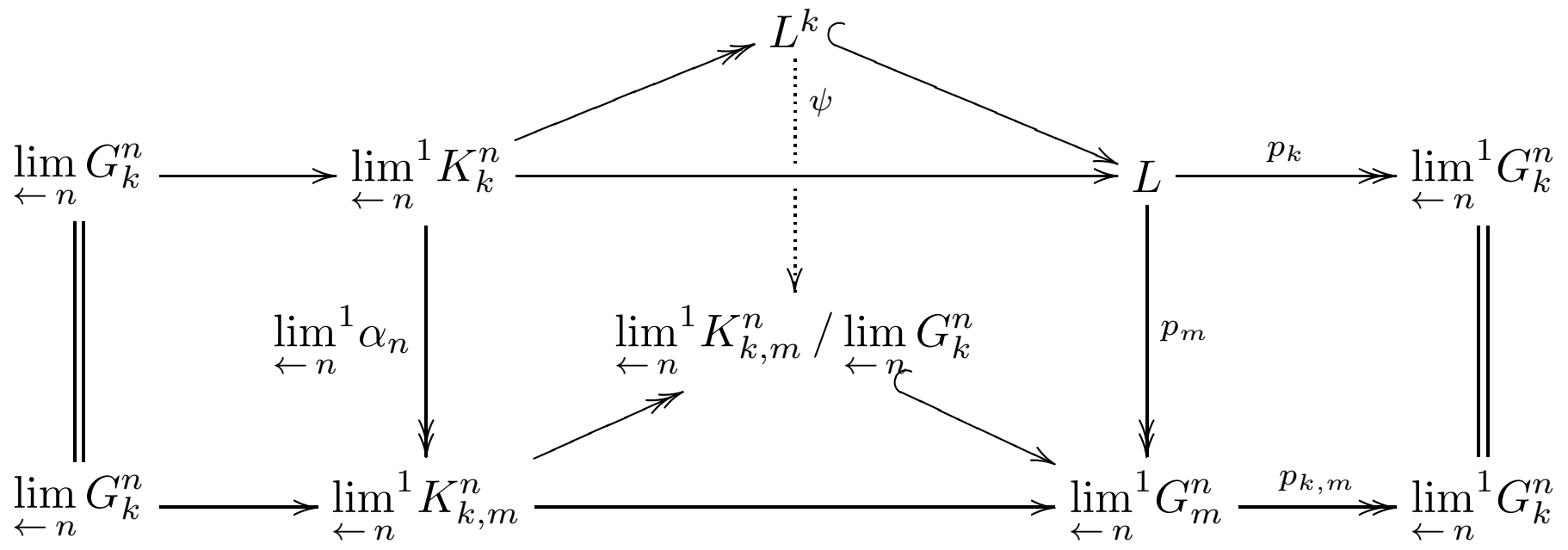

where the dotted arrow $\psi$ is to be constructed.

The upper row is the six-terms $\lim\limits_{\leftarrow}-\lim\limits_{\leftarrow}^1$ exact sequence associated to the bottom row of the *first* diagram: it identifies $L^k := \operatorname{Ker} p_k$ as the quotient of $\lim\limits_{\leftarrow n}^1 K_k^n$ by some action of the group $\lim\limits_{\leftarrow n} G_k^n$.

The bottom row is the six-terms $\lim\limits_{\leftarrow} - \lim\limits_{\leftarrow}^1$ exact sequence associated to the right column of the *first* diagram: it identifies $\operatorname{Ker} p_{k,m}$ as the quotient of $\lim\limits_{\leftarrow n}^1 K_{k,m}^n$ by some action of the group $\lim\limits_{\leftarrow n} G_k^n$.

Commutativity of the front face in the *second* diagram follows by naturality from the *first* diagram. Moreover, the map $\lim\limits_{\leftarrow n}^1 \alpha_n$ commutes with the action of $\lim\limits_{\leftarrow n} G_k^n$. It then induces the dotted map $\psi$ between the orbit sets. The map $\psi$ is easily seen to be a surjection. Also, for any $x, y \in L^k$, we have $\psi(x) = \psi(y)$ if and only if $p_m(x) = p_m(y)$. This holds precisely, by definition, if and only if $x \sim_m y$. We conclude that $\psi$ induces a bijection from $L^k/_{\sim_m}$ onto its image, and the claim is proved. ☐

*Proof of Theorem 2.15.* Suppose $L^k = L^\infty$ for some $k$. From Point *3)* of Theorem 4.4, we then deduce that $\operatorname{Ind}_k\{G_n\}_n = \infty$. By Definition 4.2 and Remark 4.3, this readily implies that $\lim\limits_{\leftarrow n}^1 K_k^n = *$. The result follows then by Point *1)* of Theorem 4.4. ☐

The rest of that section is devoted to the interpretation of Definition 4.2 in the context of $\mathcal{FD}$-phantom maps (Proposition 4.7 below).

One application is the following. As $\mathrm{Ph}^k_{\mathcal{FD}}(X, Y)$ is the image of the map $\mathrm{Ph}_{\mathcal{FD}}(X, Y\langle k\rangle) \to \mathrm{Ph}_{\mathcal{FD}}(X, Y)$, it is easy to see that if all maps in $\mathrm{Ph}_{\mathcal{FD}}(X, Y\langle k\rangle)$ have Gray index at least $m \geq k$, then $\mathrm{Ph}^k_{\mathcal{FD}}(X, Y) = \mathrm{Ph}^m_{\mathcal{FD}}(X, Y)$. The following result shows that, at least in the finite type case, the converse is true!

**Proposition 4.5.** *Let $X$ and $Y$ be nilpotent spaces of finite type, and $m \geq k$. Then the following statements are equivalent:*
*(i)* $\mathrm{Ph}^k(X, Y) = \mathrm{Ph}^m(X, Y)$*;*
*(ii) For any $i \leq k$, we have* $\mathrm{Ph}^k(X, Y\langle i\rangle) = \mathrm{Ph}^m(X, Y\langle i\rangle)$*;*
*(iii)* $\mathrm{Ph}^m(X, Y\langle k\rangle) = \mathrm{Ph}(X, Y\langle k\rangle)$.

Therefore, if $X$ and $Y$ are finite type, the Gray filtration on $\mathrm{Ph}(X, Y)$ completely determines the Gray filtration on $\mathrm{Ph}(X, Y\langle k\rangle)$, for any $k$.

**Lemma 4.6.** *Let $1 \to \{A_n\}_n \to \{B_n\}_n \to \{C_n\}_n \to 1$ be an exact sequence of towers of groups. Suppose each map $A_{n+1} \to A_n$ in the first tower is a surjection. Then* $\operatorname{Ind}\{B_n\}_n = \operatorname{Ind}\{C_n\}_n$.

*Proof.* For each $k$ we have an exact sequence of towers $1 \to \{A_k^n\}_n \to \{B_k^n\}_n \to \{C_k^n\}_n \to 1$. Indeed, if $n < k$, this is the original exact sequence $1 \to A_n \to B_n \to C_n \to 1$. If $n \geq k$, we have $A_k^n = A_k$, and the sequence $1 \to A_k \to B_k^n \to C_k^n \to 1$ is easily seen to be exact.

We apply the $\varprojlim_n^1$. As almost all maps in the tower $\{A_k^n\}_n$ are identity, we have $\varprojlim_n^1 A_k^n = *$. We deduce that the map $\varprojlim_n^1 B_k^n \to \varprojlim_n^1 C_k^n$ is a surjection with trivial kernel. Therefore its domain is $*$ if and only if its target is $*$. The result follows by definition of the index. $\square$

**Proposition 4.7.** *Let $X$ and $Y$ be nilpotent of finite type, and consider the tower $\{G_n\}_n$ where $G_n := [X, \Omega Y^{(n)}]$.*
*1) The index $\mathrm{Ind}\{G_n\}_n$ is the minimum value of $G(f)$ for $f \in \mathrm{Ph}(X,Y)$.*
*2) For any $k \geq 1$, the $k$-th index $\mathrm{Ind}_k\{G_n\}_n$ is the minimum value of $G(f)$ for $f \in \mathrm{Ph}(X, Y\langle k\rangle)$.*

*Proof.* Point *1)* follows directly from Point *3)* of Theorem 4.4: indeed it shows in particular that $\mathrm{Ph}(X,Y) = \mathrm{Ph}^k(X,Y)$ if and only if $k \leq \mathrm{Ind}\{G_n\}_n$.

Point *2)*. Using Point *1)*, it suffices to prove that the $k$-th index of $\{G_n\}_n$ is the index $\mathrm{Ind}\big\{[X, \Omega Y\langle k\rangle^{(n)}]\big\}_n$.

From the fibration sequence $\Omega^2 Y^{(k)} \to \Omega Y\langle k\rangle^{(n)} \to \Omega Y^{(n)} \to \Omega Y^{(k)}$, we deduce a surjection from $[X, \Omega Y\langle k\rangle^{(n)}]$ onto $K_k^n := \mathrm{Ker}\big([X,\Omega Y^{(n)}] \to [X, \Omega Y^{(k)}]\big)$. Let $H_k^n$ be the kernel of that surjection. Using the fibration sequence once more, we deduce a surjection from $[X, \Omega^2 Y^{(k)}]$ onto $H_k^n$.

Consider now the following diagram:

$$\begin{array}{ccccccccc}
 & 1 & \twoheadrightarrow & H_k^{n+1} & \longrightarrow & [X, \Omega Y\langle k\rangle^{(n+1)}] & \longrightarrow & K_k^{n+1} & \longrightarrow 1 \\
{[X, \Omega^2 Y^{(k)}]} \nearrow\!\!\!\!\nearrow \searrow\!\!\!\!\searrow & & & \downarrow & & \downarrow & & \downarrow & \\
 & 1 & \twoheadrightarrow & H_k^{n} & \longrightarrow & [X, \Omega Y\langle k\rangle^{(n)}] & \longrightarrow & K_k^{n} & \longrightarrow 1
\end{array}$$

where the rows are exact by definition. We see that for any $n$ the map $H_k^{n+1} \to H_k^n$ is in fact a surjection. We then apply Lemma 4.6 and deduce that $\mathrm{Ind}\{K_k^n\}_n = \mathrm{Ind}\big\{[X, \Omega Y\langle k\rangle^{(n)}]\big\}_n$. By definition, $\mathrm{Ind}\{K_k^n\}_n$ is the $k$-th index of $\{G_n\}_n$, and the result follows. $\square$

*Proof of Proposition 4.5.* Consider the following natural commutative diagram:

$$\begin{array}{ccccc}
\mathrm{Ph}^m(X,Y) & \twoheadleftarrow & \mathrm{Ph}^m(X,Y\langle i\rangle) & \twoheadleftarrow & \mathrm{Ph}^m(X,Y\langle k\rangle) \\
\downarrow & & \downarrow & & \downarrow \\
\mathrm{Ph}^k(X,Y) & \twoheadleftarrow & \mathrm{Ph}^k(X,Y\langle i\rangle) & \twoheadleftarrow & \mathrm{Ph}(X,Y\langle k\rangle)
\end{array}$$

Statement *(i)* (resp. *(ii)*, *(iii)*) means that the left (resp. middle, right) inclusion is in fact a bijection. It's then clear that *(iii)* implies *(ii)*, which itself implies *(i)*.

Suppose now that *(i)* is true. By Point *3)* of Theorem 4.4, that implies that $m \leq \mathrm{Ind}_k\{G_n\}_n$. By Proposition 4.7, we deduce that all maps $f \in \mathrm{Ph}(X, Y\langle k\rangle)$ have Gray index $G(f) \geq m$, that is that *(iii)* is true. $\square$

## 5. The completion approach to $\mathcal{F}$-phantom maps

This section is devoted to proofs of Section 2 results, excepted Point *1)* of Theorem 2.4, and Theorem 2.15, to be found in Section 4.

**Lemma 5.1.** *For any spaces $X$ and $Y$, we have $\mathrm{Ph}_{\mathcal{F}}(X,Y) = \mathrm{Ph}^1_{\mathcal{F}}(X,Y)$. Moreover, for any $k$, the 1-connected cover $Y\langle 1\rangle \to Y$ induces a bijection $\mathrm{Ph}^k_{\mathcal{F}}(X, Y\langle 1\rangle) \cong \mathrm{Ph}^k_{\mathcal{F}}(X,Y)$.*

*Proof.* This is principally Zabrodsky observation [26] that the map $\mathrm{Ph}_{\mathcal{F}}(X, Y\langle 1\rangle) \to \mathrm{Ph}_{\mathcal{F}}(X,Y)$ induced by the 1-connected cover is a bijection. It's easy to see that it induces also the stated bijections for any $k$. □

For any nilpotent finite type space $Y$, let $Y_\rho$ be the homotopy fiber of the completion of $Y$. We then have a fibration:

$$\Omega\widehat{Y} \xrightarrow{\delta} Y_\rho \xrightarrow{j} Y \xrightarrow{\widehat{e}} \widehat{Y}$$

It follows from Sullivan characterization that a map $f : X \to Y$ is a $\mathcal{F}$-phantom map if and only if it factorizes through $Y_\rho$. From Definition 2.1, it's easy to see that the Gray index $G(f) \geq k$ if and only if $f$ factorizes through the composite

$$(Y\langle k\rangle)_\rho \longrightarrow Y\langle k\rangle \longrightarrow Y$$

The following fundamental result of Roitberg-Touhey allows us a easy description of the spaces $(Y\langle k\rangle)_\rho$. Here $\widehat{\mathbb{Z}}$ denotes the profinite completion of the ring $\mathbb{Z}$ of integers.

**Theorem 5.2.** [21] *Let $Y$ be a nilpotent space of finite type, with torsion fundamental group. Then the space $Y_\rho$ splits as an infinite product of Eilenberg-MacLane spaces: $Y_\rho \simeq \Pi_{m\geq 1} K(\pi_{m+1}(Y) \otimes \widehat{\mathbb{Z}}/\mathbb{Z}, m)$.*

**Lemma 5.3.** *Let $X$ be any 0-connected space and $Y$ nilpotent of finite type. Then $H^k(X; \pi_{k+1}(Y) \otimes \widehat{\mathbb{Z}}/\mathbb{Z}) \cong 0$ if and only if $H^k(X; \pi_{k+1}(Y) \otimes \mathbb{Q}) \cong 0$. Otherwise, $H^k(X; \pi_{k+1}(Y) \otimes \widehat{\mathbb{Z}}/\mathbb{Z})$ is uncountable.*

*Proof.* Recall that $\widehat{\mathbb{Z}}/\mathbb{Z}$ is an uncountable rational group. The first assertion then follows from the obvious isomorphisms

$$\begin{aligned} H^k(X; \pi_{k+1}(Y) \otimes \widehat{\mathbb{Z}}/\mathbb{Z}) &\cong H^k(X;\mathbb{Q}) \otimes \big(\pi_{k+1}(Y) \otimes \widehat{\mathbb{Z}}/\mathbb{Z}\big) \\ &\cong \big(H^k(X;\mathbb{Q}) \otimes \pi_{k+1}(Y) \otimes \mathbb{Q}\big) \otimes \widehat{\mathbb{Z}}/\mathbb{Z} \cong H^k(X; \pi_{k+1}(Y) \otimes \mathbb{Q}) \otimes \widehat{\mathbb{Z}}/\mathbb{Z}. \end{aligned}$$

The choice of any non zero element of $H^k(X; \pi_{k+1}(Y) \otimes \mathbb{Q})$ defines an inclusion $\mathbb{Q} \hookrightarrow H^k(X; \pi_{k+1}(Y) \otimes \mathbb{Q})$ which, once tensorized by $\widehat{\mathbb{Z}}/\mathbb{Z}$, gives an inclusion $\widehat{\mathbb{Z}}/\mathbb{Z} \hookrightarrow H^k(X; \pi_{k+1}(Y) \otimes \widehat{\mathbb{Z}}/\mathbb{Z})$. The second assertion follows. □

*Proof of Points 2) and 3) in Theorem 2.4.* By Lemma 5.1 we can always suppose $k \geq 1$. Consider the following diagram:

$$\begin{array}{ccccc} & & & & (Y\langle k+1\rangle)_\rho \\ & & & \swarrow & \big\downarrow \\ X & \xrightarrow{f} & Y & & \\ & & & \nwarrow & \\ & & & & (Y\langle k\rangle)_\rho \simeq (Y\langle k+1\rangle)_\rho \times K(\pi_{k+1}(Y) \otimes \widehat{\mathbb{Z}}/\mathbb{Z}, k) \end{array}$$

where the vertical map is merely the inclusion of the first factor. Here we use Theorem 5.2 to describe the spaces $(Y\langle k\rangle)_\rho$.

If $H^k(X, \pi_{k+1}(Y) \otimes \widehat{\mathbb{Z}}/\mathbb{Z}) \cong 0$, then a map $f$ factors through $(Y\langle k\rangle)_\rho$ if and only if it factors through $(Y\langle k+1\rangle)_\rho$, and $\mathrm{Ph}^k_{\mathcal{F}}(X,Y) = \mathrm{Ph}^{k+1}_{\mathcal{F}}(X,Y)$.

If $H^m(X, \pi_{m+1}(Y) \otimes \widehat{\mathbb{Z}}/\mathbb{Z}) \cong 0$ for any $m \geq k$, then the set $[X, (Y\langle k\rangle)_\rho] \cong \Pi_{m\geq k}\, [X, K(\pi_{m+1}(Y) \otimes \widehat{\mathbb{Z}}/\mathbb{Z}, m)]$ is the trivial set, and $\mathrm{Ph}^k_{\mathcal{F}}(X,Y) = \{*\}$. □

The proof of Theorem 2.7 begins with the following lemma:

**Lemma 5.4.** *Let $f \in \mathrm{Ph}_{\mathcal{F}}(X,Y)$, where $Y$ is a nilpotent space of finite type, with torsion fundamental group, and $X$ is any space. Then the following statements are equivalent:*
*(i) $G(f) \geq k$;*
*(ii) There exists $\phi : X \to Y_\rho$ a lift of $f$ such that $X \xrightarrow{\phi} Y_\rho \to (Y_\rho)^{(k-1)}$ is inessential;*
*(iii) If $\phi : X \to Y_\rho$ is* any *lift of $f$, then there exists a map $\lambda$ making commutative the following diagram:*

$$\begin{array}{ccccc} X & \xrightarrow{\phi} & Y_\rho & \xrightarrow{j} & Y \\ {\scriptstyle\lambda}\downarrow\vdots & & \downarrow{\scriptstyle\pi} & & \\ \Omega\widehat{Y} & \xrightarrow{\delta} & Y_\rho & \xrightarrow{\pi} & (Y_\rho)^{(k-1)} \end{array} \qquad (f = j\circ\phi : X \to Y)$$

*Proof.* As Theorem 5.2 applies, we have easy identifications $(Y^{(k)})_\rho \simeq (Y_\rho)^{(k-1)}$ and $(Y\langle k\rangle)_\rho \simeq (Y_\rho)\langle k-1\rangle$, and a (trivial) fibration $(Y\langle k\rangle)_\rho \xrightarrow{i} Y_\rho \xrightarrow{\pi} (Y_\rho)^{(k-1)}$.

Moreover, we have a commutative diagram:

$$\begin{array}{ccccccc} \Omega\widehat{Y} & \xrightarrow{\delta} & Y_\rho & \xrightarrow{j} & Y & \xrightarrow{\widehat{e}} & \widehat{Y} \\ {\scriptstyle\pi'}\downarrow & & {\scriptstyle\pi}\downarrow & & \downarrow & & \downarrow \\ (\Omega\widehat{Y})^{(k-1)} \simeq \Omega\widehat{Y}^{(k)} & \xrightarrow{\delta^{(k-1)}} & (Y_\rho)^{(k-1)} \simeq (Y^{(k)})_\rho & \longrightarrow & Y^{(k)} & \longrightarrow & \widehat{Y}^{(k)} \end{array}$$

where the rows are fibrations and all vertical maps are Postnikov sections. By naturality [8, Proposition 11.5], we also have a commutative diagram linking the actions $A$ and $A^{k-1}$ induced by these fibrations, namely, for any space $X$:

$$\begin{array}{ccc} [X,\Omega\widehat{Y}] \times [X,Y_\rho] & \xrightarrow{A} & [X,Y_\rho] \\ {\scriptstyle\pi'_*\times\pi_*}\downarrow & & \downarrow{\scriptstyle\pi_*} \\ [X,(\Omega\widehat{Y})^{(k-1)}] \times [X,(Y_\rho)^{(k-1)}] & \xrightarrow{A^{k-1}} & [X,(Y_\rho)^{(k-1)}] \end{array}$$

Suppose $G(f) \geq k$. Choose $\psi : X \to (Y\langle k\rangle)_\rho$ any lift of $f$. Then $i\circ\psi : X \to Y_\rho$ is a lift of $f$ projecting to $*$ in $(Y_\rho)^{(k-1)}$, and *(ii)* follows. Let now $\phi : X \to Y_\rho$ be a lift of $f$. We have $j\circ i\circ\psi \simeq f \simeq j\circ\phi$. Then, by classical theory, there exists $\lambda : X \to \Omega\widehat{Y}$ such that $\phi = A(\lambda, i\circ\psi)$. We deduce $\pi\circ\phi \simeq A^{k-1}(\pi'\circ\lambda, \pi\circ i\circ\psi) \simeq A^{k-1}(\pi'\circ\lambda, *) \simeq \delta^{(k-1)}\circ\pi'\circ\lambda \simeq \pi\circ\delta\circ\lambda$, and *(iii)* follows.

Suppose *(iii)* is true. Choose $\phi$ and $\lambda$ such that $\pi\circ\phi \simeq \pi\circ\delta\circ\lambda$. Set $\tilde{\phi} := A(-\lambda,\phi) : X \to Y_\rho$. Then $\tilde{\phi}$ is also a lift of $f$. We have $\pi\circ\tilde{\phi} \simeq A^{k-1}(\pi'\circ(-\lambda), \pi\circ\phi) \simeq A^{k-1}(-(\pi'\circ\lambda), \delta^{(k-1)}\circ\pi'\circ\lambda) \simeq \delta^{(k-1)}\circ(-(\pi'\circ\lambda)+\pi'\circ\lambda) \simeq *$ and *(ii)* follows (the signs $+$ and $-$ above refer to the group structure on $[X,\Omega\widehat{Y}]$ and $[X,(\Omega\widehat{Y})^{(k-1)}]$). This shows moreover that there exists $\psi : X \to (Y\langle k\rangle)_\rho$ such that $i\circ\psi = \tilde{\phi}$. In particular $\psi$ is a lift of $f$, and $G(f) \geq k$. □

*Proof of Theorem 2.7.* By hypothesis, the domain $X$ is 0-connected, and then we have bijections $[X,\Omega\widehat{Y}\langle 1\rangle] \cong [\Sigma X,\widehat{Y}\langle 1\rangle] \cong [\Sigma X,\widehat{Y}] \cong [X,\Omega\widehat{Y}]$. Together with Lemma 5.1, this implies we can always suppose that $Y$ is 1-connected. We can then use Theorem 5.2 to describe $Y_\rho$.

Point *1)*. Let $f \in \mathrm{Ph}^{\infty}_{\mathcal{F}}(X,Y)$, and choose $\phi : X \to Y_\rho$ any lift of $f$. By Theorem 5.2, we can identify $\phi$ with the collection $(\phi^k)_k$, where $\phi^k$ is the composition $X \xrightarrow{\phi} Y_\rho \to (Y_\rho)^{(k)}$.

Let now $\lambda_1, \dots, \lambda_s$, be the finitely many elements of $[X, \Omega\widehat{Y}]$. Similarly, for any $1 \le m \le s$, we identify $\delta \circ \lambda_m : X \to Y_\rho$ with the collection $(\lambda^k_m)_k$, where $\lambda^k_m$ is the composition $X \xrightarrow{\lambda_m} \Omega\widehat{Y} \xrightarrow{\delta} Y_\rho \longrightarrow (Y_\rho)^{(k)}$.

Suppose $f$ is essential. Then for any $1 \le m \le s$, the maps $\phi$ and $\delta \circ \lambda_m$ are non homotopic. Then there exists some integer $k(m)$, depending on $m$, such that $\phi^{k(m)}$ and $\lambda^{k(m)}_m$ are non homotopic. Take now any integer $k-1$ greater than all the $k(m)$. Then for any $1 \le m \le s$, the maps $\phi^{k-1}$ and $\lambda^{k-1}_m$ are non homotopic. By Lemma 5.4, we deduce $G(f) < k$, a contradiction. Therefore $f$ is inessential and the result follows.

For a distinct proof of the same result, see Remark 5.12.

Point *2)*. Suppose $H^k(X; \pi_{k+1}(Y) \otimes \mathbb{Q}) \not\cong 0$. To show that the inclusion $\mathrm{Ph}^k_{\mathcal{F}}(X,Y) \supset \mathrm{Ph}^{k+1}_{\mathcal{F}}(X,Y)$ is strict, we will produce a $\mathcal{F}$-phantom map $f : X \to Y$ with Gray index exactly $k$.

Let $K := K(\pi_{k+1}(Y) \otimes \widehat{\mathbb{Z}}/\mathbb{Z}, k)$. By Lemma 5.3, we see that $[X, K]$ is uncountable. By hypothesis, $[X, \Omega\widehat{Y}]$ is countable, and we can choose a map $\phi : X \to K$ which *does not* factor through $\Omega\widehat{Y} \to Y_\rho \xrightarrow{proj.} K$. Let $f$ be the composition $X \xrightarrow{\phi} K \hookrightarrow Y_\rho \to Y$.

The composition $X \xrightarrow{\phi} K \hookrightarrow Y_\rho \xrightarrow{\pi} (Y_\rho)^{(k-1)}$ is clearly trivial, and by Lemma 5.4 we deduce that $G(f) \ge k$.

Suppose now that $G(f) \ge k+1$. By Lemma 5.4, there exists then a map $\lambda$ making commutative the following diagram:

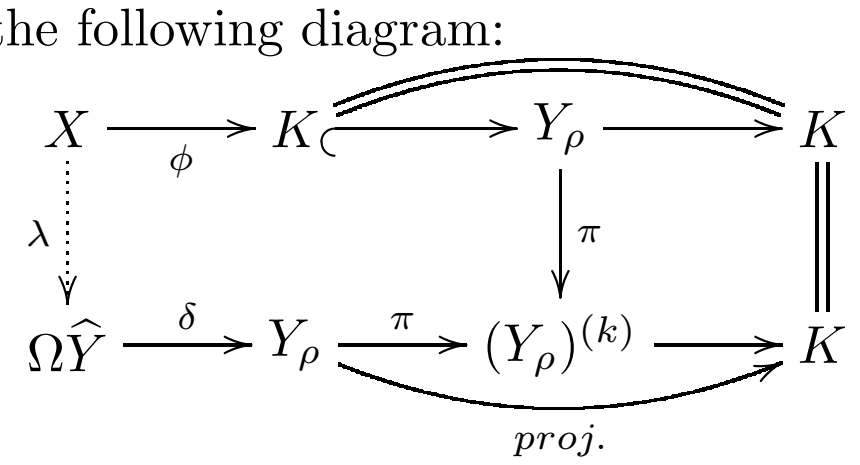

This clearly contradicts our choice of the map $\phi$, and then $G(f) = k$. ☐

We list below some already known cases where Theorem 2.7 do apply.

**Proposition 5.5.** *In each of the following cases, we have* $[X, \Omega\widehat{Y}] \cong *$*:*
*1) The spaces* $X$ *and* $Y$ *are finite type, with torsion fundamental group, and* $X \simeq \Sigma^i Z^{(m)}$ *whereas* $Y \simeq \Omega^j W$*, for some integers* $m$ *and* $i, j \ge 0$*, for some space* $Z$ *and some finite CW-complex* $W$*.*
*2) The space* $\Omega Y$ *is a one-connected finite CW-complex, and* $X \simeq BG$*, the classifying space of some 0-connected Lie group* $G$ [12, Theorem 5.6.(ii)].
*3) The space* $\Omega Y$ *is a finite CW-complex, and* $X$ *is a 0-connected infinite loop space with torsion fundamental group* [13, Theorem 3].

*Proof.* As far as we know, Point *1)* is not clearly stated in the litterature, although undoubtely known by experts. Here is a quick proof. Let $X_\tau$ be the homotopy fiber of the rationalization $X \to X_{(0)}$. Then $[\Sigma X_\tau, Y] \cong *$, by [26, Theorem D]. But $[\Sigma X_\tau, Y] \cong [X, \Omega\widehat{Y}]$ by [23, Proposition 34], and Point *1)* follows. ☐

**Remark 5.6.** In the particular case where $[X, \Omega\widehat{Y}]$ is trivial, we have an easy description of the Gray filtration on $\mathrm{Ph}_{\mathcal{F}}(X,Y)$. Indeed, by Theorem 5.2, we have

a bijection:

$$\mathrm{Ph}_{\mathcal{F}}(X,Y) \cong [X, Y_\rho] \cong \Pi_{m\geq 1} H^m(X; \pi_{m+1}(Y)\otimes \widehat{\mathbb{Z}}/\mathbb{Z}).$$

Through this product description, we see thanks to Lemma 5.4 that a phantom map has Gray index $k$ if and only if its $(k-1)$-th first coordinates are zero.

We shall now go towards the proof of Theorem 2.8. Recall that we have defined spaces $P^nY$ as homotopy pullbacks:

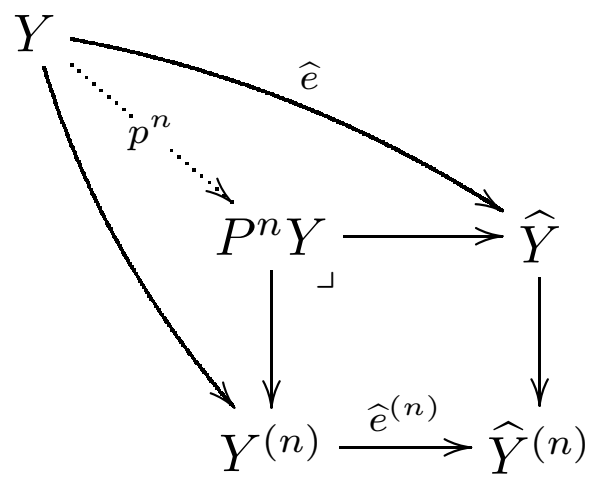

There might exist several homotopy classes of maps $p^n$ making the above diagram commutative. To avoid this ambiguity, we shall state the following lemma.

**Lemma 5.7.** *For any $n \geq 0$, there exists a functor $P^n$ from spaces to spaces, together with a coaugmentation (a natural transformation) $p^n : Id \to P^n$, such that, for any space $Y$, the space $P^nY$ is the homotopy pullback associated to the maps $Y^{(n)} \xrightarrow{\widehat{e}^{(n)}} \widehat{Y}^{(n)} \longleftarrow \widehat{Y}$.*

*Proof.* Let $F$ be a functor from spaces to spaces with coaugmentation $\eta : Id \to F$. Such a functor associates, to any map $f : X \to Y$, a square

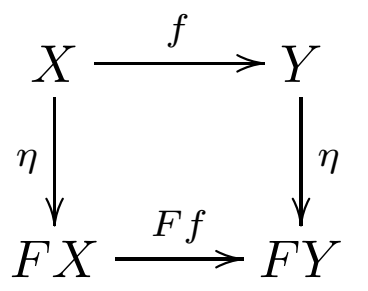

which is *strictly commutative*, and not just commutative up to homotopy.

The technical key point of our proof is the fact that Postnikov section $p^{(n)} : Y \to Y^{(n)}$, as well as completion $\widehat{e} : Y \to \widehat{Y}$, can be chosen as such functors. This comes from the theory of localization with respect to a map, see for example the lines following Theorem 6.3. of [2].

Having chosen these functors (Postnikov section and completion), we then define $P^nY$, as usual, as the set of triples $(x, \omega, z)$ where $x \in Y^{(n)}$, $z \in \widehat{Y}$, and $\omega$ is a path in $\widehat{Y}^{(n)}$, starting at $\widehat{e}^{(n)}(x)$, ending at $p^{(n)}(z)$. We also define $p^n : Y \to P^nY$ by $y \mapsto \big(p^{(n)}(y), c_y, \widehat{e}(y)\big)$, where $c_y$ denotes the constant path at the point $\widehat{e}^{(n)}p^{(n)}(y) = p^{(n)}\widehat{e}(y)$. For any map $f : X \to Y$, strict commutativity shows that the assignement

$$\big(x, t \mapsto \omega(t), z\big) \mapsto \big(f^{(n)}(y), t \mapsto \widehat{f}^{(n)}(\omega(t)), \widehat{f}(z)\big)$$

well defines a map $P^nf : P^nX \to P^nY$ such that the following square strictly commutes:

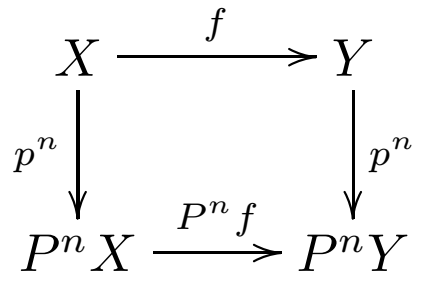

□

**Remark 5.8.** More generally, to any pair of maps $\alpha : A \to B$ and $\beta : C \to D$, we can associate a functor $L_{\alpha,\beta}$ from spaces to spaces, together with a coaugmentation $\eta_{\alpha,\beta} : Id \to L_{\alpha,\beta}$, as follows. Let $L_\alpha$ and $L_\beta$ be the localization functors with respect to the maps $\alpha, \beta$, together with their coaugmentation $\eta_\alpha, \eta_\beta$. Then the functor $L_{\alpha,\beta}$ takes a space $Y$ into the homotopy pullback associated to the diagram:

$$L_\alpha Y \overset{L_\alpha(\eta_\beta)}{\longrightarrow} L_\alpha(L_\beta Y) \overset{\eta_\alpha}{\longleftarrow} L_\beta Y$$

In the proof above $\alpha$ is Postnikov section and $\beta$ is completion.

**Lemma 5.9.** *1) The map $p^n_* : \pi_k(Y) \to \pi_k(P^n Y)$ induced in homotopy by $p^n$ is an isomorphism if $k \le n$. If $k > n$, it is, up to natural isomorphism of its target, the completion $\widehat{e}_* : \pi_k(Y) \to \pi_k(\widehat{Y})$.*
*2) For any space $Y$ and $k \ge n$, the map $p^k : P^n Y \to P^k(P^n Y)$ is an equivalence.*
*3) For any map $f : X \to P^n Y$ and $k \ge n$, there exists, up to homotopy, a natural factorization of $f$ as shown in the following diagram:*

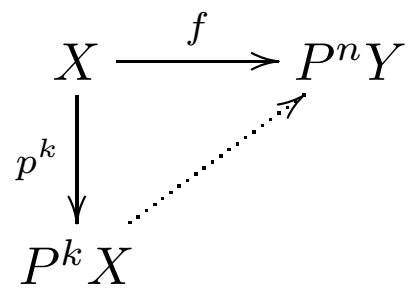

*Proof.* Point *1)* follows from the Mayer-Vietoris sequence in homotopy of the pullback involving $P^n Y$, and implies easily Point *2)*. In the following diagram,

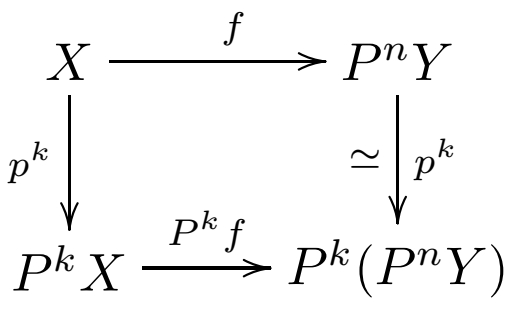

the right vertical map $p^k$ is then an equivalence when $k \ge n$. If $(p^k)^{-1}$ denotes its homotopy inverse, the composition $(p^k)^{-1} \circ P^k f$ is the desired factorization of $f$. ☐

**Corollary 5.10.** *There exists a commutative diagram*

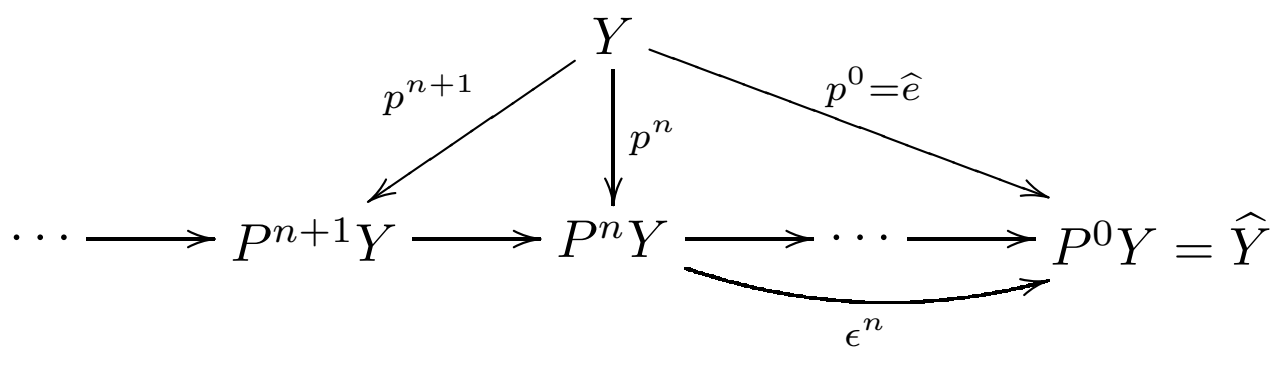

*and the maps $p^n$ induce a equivalence from $Y$ to the homotopy inverse limit of the tower $\{P^n Y\}_n$.*

*Proof of Theorem 2.8.* Construct the following diagram, where all vertical sequences are fibrations:

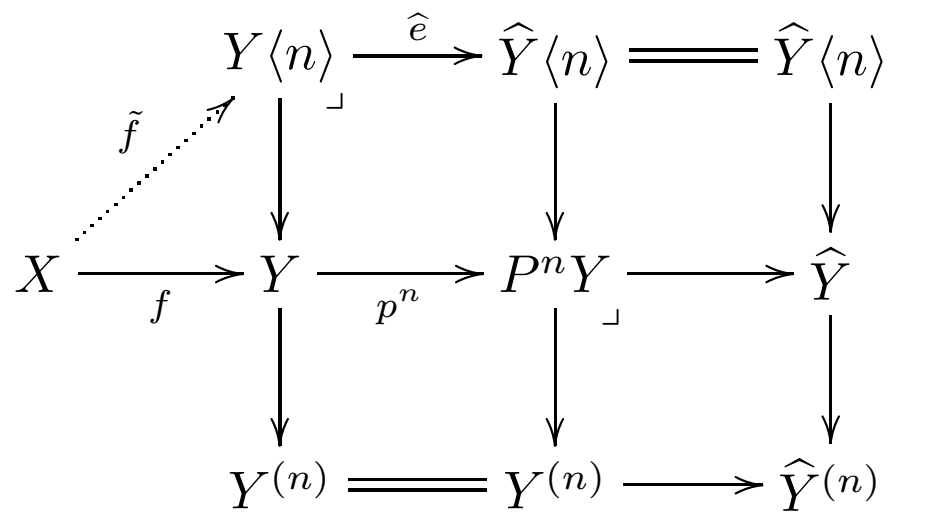

The left upper square is then a pullback, by construction.

Let $f \in \mathrm{Ph}_{\mathcal{F}}(X,Y)$. If there exists a lift $\tilde{f}$ such that $\widehat{e} \circ \tilde{f} \simeq *$, then certainly $p^n \circ f \simeq *$. But, by pullback property, the converse is also true. As $\widehat{e} \circ \tilde{f} \simeq *$ if and only if $\tilde{f} \in \mathrm{Ph}_{\mathcal{F}}(X, Y\langle n\rangle)$, Point *1)* follows.

We deduce the existence of an exact sequence of pointed sets

$$* \longrightarrow \mathrm{Ph}^{\infty}_{\mathcal{F}}(X,Y) \longrightarrow [X,Y] \longrightarrow \varprojlim_{n} [X, P^nY] \longrightarrow *$$

the base point being the trivial map. Then Point *2)* follows from [1, Ch. IX, Corollary 3.2.]. □

**Corollary 5.11.** *Let $Y$ be a nilpotent space of finite type, and $X$ be any space. Set*

$$H_n := \mathrm{Im}\left([X, \Omega P^nY] \to [X, \Omega Y^{(n)}]\right) \quad \textit{and} \quad \widehat{H}_n := \mathrm{Im}\left([X, \Omega P^nY] \to [X, \Omega \widehat{Y}]\right)$$

*1) The surjections $[X, \Omega P^nY] \twoheadrightarrow H_n$ induce a bijection $\mathrm{Ph}^{\infty}_{\mathcal{F}}(X,Y) \cong \varprojlim_{n}^1 H_n$.*

*2) The surjections $[X, \Omega P^nY] \twoheadrightarrow \widehat{H}_n$ induce a bijection $\mathrm{Ph}^{\infty}_{\mathcal{F}}(X,Y) \cong \varprojlim_{n}^1 \widehat{H}_n$.*

**Remark 5.12.** Point *2)* of Corollary 5.11 gives an other proof of Point *1)* of Theorem 2.7. Indeed, if $[X, \Omega\widehat{Y}]$ is finite, so are its subgroups $\widehat{H}_n$. We deduce $\varprojlim_{n}^1 \widehat{H}_n \cong *$, and the result follows.

*Proof of Corollary 5.11.* From Lemma 5.1, we can check that we can always suppose that $Y$ is one-connected.

Point *1)*. We consider the diagram:

$$\begin{array}{ccc} [X,\Omega P^nY] & \twoheadrightarrow & H_n \\ & \searrow & \downarrow \\ & & [X, \Omega Y^{(n)}] \end{array} \qquad\qquad \begin{array}{ccc} \mathrm{Ph}^{\infty}_{\mathcal{F}}(X,Y) & \twoheadrightarrow & \varprojlim_{n}^1 H_n \\ & \searrow & \downarrow \\ & & \mathrm{Ph}_{\mathcal{F}D}(X,Y) \end{array}$$

whose right part is obtained from the left one, by taking $\varprojlim_{n}^1$. As the map $\mathrm{Ph}^{\infty}_{\mathcal{F}}(X,Y) \to \mathrm{Ph}_{\mathcal{F}D}(X,Y)$ is merely an inclusion between two subsets of $[X,Y]$, the assertion follows.

Point *2)*. Let $A_n$ be the image of $[X, \Omega\widehat{Y}] \to [X, (Y_\rho)^{(n-1)}]$ (notice that in general $A_n$ as no reason to be a group). We deduce from Lemma 5.4 that the image of the composition

$$\varprojlim_{n} A_n \hookrightarrow \varprojlim_{n} [X, (Y_\rho)^{(n)}] \cong [X, Y_\rho] \longrightarrow [X,Y]$$

is exactly $\mathrm{Ph}^{\infty}_{\mathcal{F}}(X,Y)$. This set then appears as an orbit set,

$$\mathrm{Ph}^{\infty}_{\mathcal{F}}(X,Y) \cong \varprojlim_{n} A_n \,/\, [X, \Omega\widehat{Y}],$$

using the action of $[X, \Omega\widehat{Y}]$ on $\varprojlim_{n} A_n$ induced by its action on $[X, Y_\rho]$.

On another hand, a carefull reading of the proof of [14, Theorem 2] reveals the following purely algebraic fact: the group $[X, \Omega\widehat{Y}]$ acts naturally on $\varprojlim_{n}\big([X, \Omega\widehat{Y}]/\widehat{H}_n\big)$ in such a way that

$$\varprojlim_{n}^1 \widehat{H}_n \cong \Big(\varprojlim_{n}\big([X, \Omega\widehat{Y}]/\widehat{H}_n\big)\Big) \,/\, [X, \Omega\widehat{Y}].$$

To conclude, it suffices to produce natural bijections $[X, \Omega\widehat{Y}]/\widehat{H}_n \to A_n$, which are moreover maps of sets on which $[X, \Omega\widehat{Y}]$ acts, that is to say maps commuting with the action. This can be done thanks to the exact sequence $1 \to \widehat{H}_n \to$ $[X,\Omega\widehat{Y}] \twoheadrightarrow A_n$ induced by the fibration sequence $\Omega P^nY \to \Omega\widehat{Y} \to (Y^{(n)})_\rho \simeq (Y_\rho)^{(n-1)}$. ☐

Theorem 2.8 as well as Corollary 5.11 give $\lim\limits_{\leftarrow n}{}^1$ descriptions of the set $\mathrm{Ph}^\infty_{\mathcal{F}}(X,Y)$. To each of these descriptions is associated an algebraic Gray filtration (Definition 2.13). The following result interprets geometrically these filtrations.

**Proposition 5.13.** *Let $Y$ be a nilpotent space of finite type, and $X$ be any space.*
*1) The algebraic Gray filtration on $\mathrm{Ph}^\infty_{\mathcal{F}}(X,Y) \cong \lim\limits_{\leftarrow n}{}^1[X,\Omega P^nY]$ is trivial.*
*2) Let $\mathrm{Ph}^{\infty,k}_{\mathcal{F}}(X,Y)$ be the $k$-th term of the algebraic Gray filtration on $\mathrm{Ph}^\infty_{\mathcal{F}}(X,Y) \cong \lim\limits_{\leftarrow n}{}^1 H_n$. Then $\mathrm{Ph}^{\infty,k}_{\mathcal{F}}(X,Y)$ is the image of the map $\mathrm{Ph}^\infty_{\mathcal{F}}(X,Y\langle k\rangle) \to \mathrm{Ph}^\infty_{\mathcal{F}}(X,Y)$ induced by the $k$-connected cover $Y\langle k\rangle \to Y$.*

*Proof.* Point *1)*. We shall prove that the algebraic Gray filtration on the $\lim\limits_{\leftarrow n}{}^1$ of the tower

$$\cdots \to [X,\Omega P^nY] \to [X,\Omega P^{n-1}Y] \to \cdots \to [X,\Omega P^1Y] \to [X,\Omega\widehat{Y}] \to *$$

is trivial. The first term $L^1$ of that filtration is by definition the kernel of the natural map $\mathrm{Ph}^\infty_{\mathcal{F}}(X,Y) \twoheadrightarrow \lim\limits_{\leftarrow n}{}^1\widehat{H}_n$. Then $L^1$ is trivial by Point *2)* of Corollary 5.11.

Point *2)*. Let $H_{n,k}$ be the image of the map $[X,\Omega P^nY\langle k\rangle] \to [X,\Omega Y\langle k\rangle^{(n)}]$. By naturality, we see that the map $[X,\Omega Y\langle k\rangle^{(n)}] \to [X,\Omega Y^{(n)}]$ induced by $Y\langle k\rangle \to Y$ takes $H_{n,k}$ into $H_n$. We claim that the sequence of groups $H_{n,k} \to H_n \to H_k$ is exact, where $H_n \to H_k$ is induced by $[X,\Omega Y^{(n)}] \to [X,\Omega Y^{(k)}]$.

Assuming that claim for the moment, consider the diagram:

$$\begin{array}{ccccc}
H_{n,k} & \longrightarrow & H_n & \longrightarrow & H_k \\
\downarrow & & \| & & \uparrow \\
K_n := \mathrm{Im}(H_{n,k} \to H_n) & \hookrightarrow & H_n & \twoheadrightarrow & H^n_k := \mathrm{Im}(H_n \to H_k)
\end{array}$$

Taking $\lim\limits_{\leftarrow n}{}^1$, we obtain

$$\begin{array}{ccccc}
\lim\limits_{\leftarrow n}{}^1 H_{n,k} \cong \mathrm{Ph}^\infty_{\mathcal{F}}(X,Y\langle k\rangle) & & & & \\
\downarrow & \searrow & & & \\
\lim\limits_{\leftarrow n}{}^1 K_n & \longrightarrow & \mathrm{Ph}^\infty_{\mathcal{F}}(X,Y) & \xrightarrow{p_k} & \lim\limits_{\leftarrow n}{}^1 H^n_k
\end{array}$$

where the row is an exact sequence of pointed sets. The $k$-th term of the filtration is by definition the kernel of $p_k$, which then clearly coincide with the image of $\mathrm{Ph}^\infty_{\mathcal{F}}(X,Y\langle k\rangle) \to \mathrm{Ph}^\infty_{\mathcal{F}}(X,Y)$.

We show now our claim. From the following diagram,

$$\begin{array}{ccccc}
H_{n,k} & \longrightarrow & H_n & \longrightarrow & H_k \\
\downarrow & & \downarrow & & \downarrow \\
[X,\Omega Y\langle k\rangle^{(n)}] & \longrightarrow & [X,\Omega Y^{(n)}] & \longrightarrow & [X,\Omega Y^{(k)}]
\end{array}$$

where the bottom row is exact, we see we have only to prove that a map $f \in H_n$ projecting to $*$ in $H_k$ belongs to the image of $H_{n,k}$.

For such an $f$, choose $g$ a lift of $f$ through $\Omega P^nY$ (a lift exists by definition of $H_n$), and consider the following diagram, where we construct the maps $u_1$, $u_2$ and

$u_3$ successively:

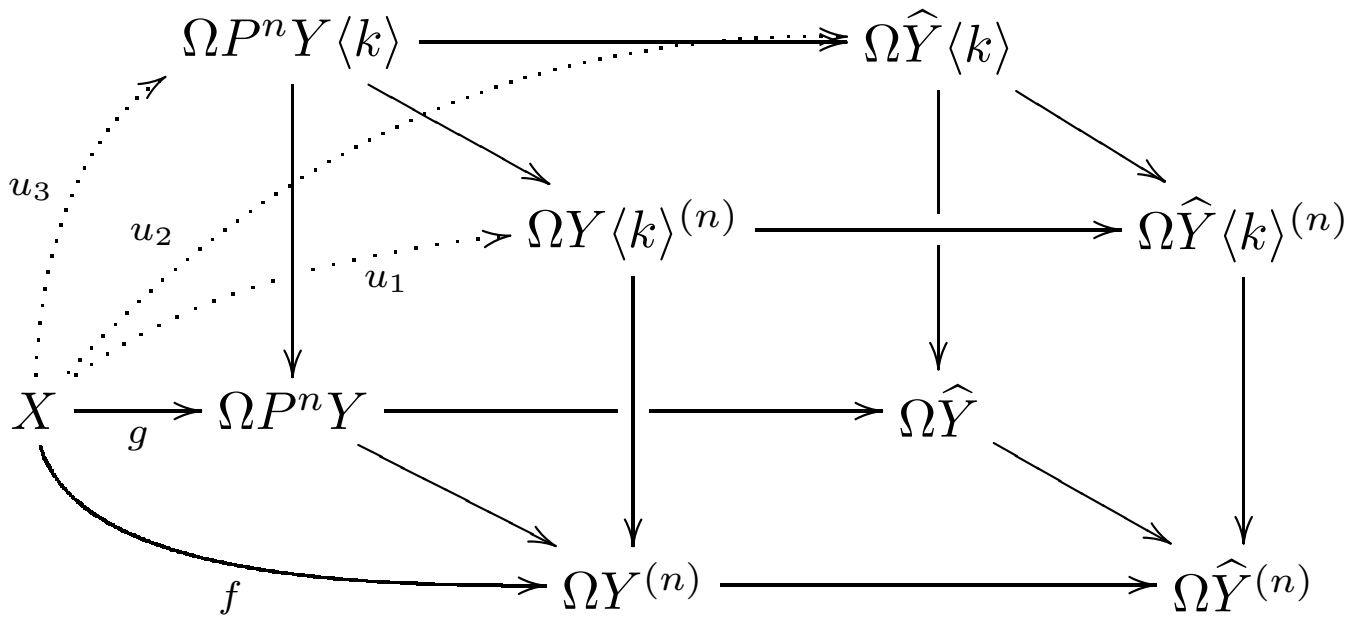

As $f$ projects to $*$ in $H_k$, there exists a lift $u_1$ of $f$ through $\Omega Y\langle k\rangle^{(n)}$. We then choose $u_2$ thanks to the pullback property of the right vertical face of the cube. Again, we choose $u_3$ thanks to the pullback property of the top horizontal face of the cube. In particular, $u_3$ is a lift of $u_1$ through $\Omega P^n Y\langle k\rangle$ (we don't claim however that $u_3$ is a lift of $g$). That means by definition that $u_1 \in H_{n.k}$, and then $f$ belongs to the image of $H_{n,k}$, as required. □

Heading now towards the proof of Theorem 2.9, we shall before make a brief review of various existing *universal phantom maps.*

Let $X = \bigcup_{n\geq 1} X_n$ be a CW-complex, where $X_n$ denotes the $n$-skeleton of $X$. The *universal $\mathcal{FD}$-phantom map out of* $X$ [6] is the right map in the cofiber sequence:

$$\bigvee_{n\geq 1} X_n \longrightarrow X \longrightarrow \bigvee_{n\geq 1} \Sigma X_n$$

It's a $\mathcal{FD}$-phantom map and, for any space $Y$, a map $f : X \to Y$ is $\mathcal{FD}$-phantom if and only if $f$ factorizes as some composition $f : X \to \bigvee_{n\geq 1} \Sigma X_n \to Y$.

Dually, for any space $Y$, the *universal $\mathcal{FD}$-phantom map into* $Y$ [11] is the left map in the fiber sequence:

$$\Pi_{n\geq 1}\, \Omega Y^{(n)} \longrightarrow Y \longrightarrow \Pi_{n\geq 1}\, Y^{(n)}$$

(by fiber sequence it is meaned here that the homotopy fiber of the right map has the weak homotopy type of $\Pi_{n\geq 1}\, \Omega Y^{(n)}$). It's a $\mathcal{FD}$-phantom map and, for any space $X$, a map $f : X \to Y$ is $\mathcal{FD}$-phantom if and only if $f$ factorizes as some composition $f : X \to \Pi_{n\geq 1}\, \Omega Y^{(n)} \to Y$.

Recall now from Section 2 the canonical map $\Phi_X : X \to W(X)$ of McGibbon-Strom [18], associated to any space $X$. It's the *universal $\mathcal{FD}$-phantom map with infinite Gray index out of* $X$. That is that $\Phi_X \in \mathrm{Ph}^{\infty}_{\mathcal{FD}}(X, W(X))$ and, for any space $Y$, a map $f : X \to Y$ belongs to $\mathrm{Ph}^{\infty}_{\mathcal{FD}}(X, Y)$ if and only if $f$ factorizes as some composition $f : X \to W(X) \to Y$.

We now complete this list, defining below two other universal phantom maps. For any space $Y$ and $k \geq 1$, set $Z_k(Y) := \Pi_{n\geq 1}\, \Omega Y\langle k\rangle^{(n)}$, and let $\theta_k(Y)$ be the composition $Z_k(Y) \to Y\langle k\rangle \to Y$, where the first map is the universal $\mathcal{FD}$-phantom map into $Y\langle k\rangle$.

**Proposition 5.14. and Definition.**
*1) Let $Y$ be any space. The* universal $\mathcal{FD}$-phantom map with infinite Gray index

into $Y$ *is the top map* $\theta_Y$ *in the following pullback square:*

$$\begin{array}{ccc} Z(Y) & \xrightarrow{\theta_Y} & Y \\ \downarrow & \lrcorner & \downarrow{\scriptstyle \Delta} \\ \Pi_{k\geq 1} Z_k(Y) & \xrightarrow[\Pi_{k\geq 1}\theta_k(Y)]{} & \Pi_{k\geq 1} Y \end{array}$$

*Then* $\theta_Y \in \mathrm{Ph}^{\infty}_{\mathcal{F}D}(Z(Y), Y)$ *and, for any space* $X$*, a map* $f : X \to Y$ *belongs to* $\mathrm{Ph}^{\infty}_{\mathcal{F}D}(X, Y)$ *if and only if* $f$ *factorizes as some composition* $f : X \to Z(Y) \to Y$*.*

*2) Let* $Y$ *be a nilpotent space of finite type. The* universal $\mathcal{F}$-phantom map with infinite Gray index into $Y$ *is the left map in the fiber sequence*

$$\Pi_{n\geq 1} \Omega P^n Y \longrightarrow Y \longrightarrow \Pi_{n\geq 1} P^n Y,$$

*where the right map is defined by the collection* $\{p^n : Y \to P^n Y\}_n$*. This map is a* $\mathcal{F}$*-phantom map with infinite Gray index and, for any space* $X$*, a map* $f : X \to Y$ *belongs to* $\mathrm{Ph}^{\infty}_{\mathcal{F}}(X, Y)$ *if and only if* $f$ *factorizes as some composition* $f : X \to \Pi_{n\geq 1} \Omega P^n Y \to Y$*.*

*Proof.* We don't give the details for Point *1)*: this is a straightforward dualization of McGibbon-Strom construction [18].

To prove that the homotopy fiber of $Y \to \Pi_{n\geq 1} P^n Y$ is weakly equivalent to $\Pi_{n\geq 1} \Omega P^n Y$, we follow the same line of arguments than in [11, Theorem 3]. The universal property is easily checked with Point *1)* of Theorem 2.8. $\square$

*Proof of Theorem 2.9.* Consider the following diagram, where the bottom fibration sequence is obtained by CW-approximation of the upper fibration sequence:

$$\begin{array}{ccccccccc} \Omega Y & \xrightarrow{\Omega\{p^n\}_n} & \Pi_{n\geq 1} \Omega P^n Y & \longrightarrow & \Pi_{n\geq 1} \Omega P^n Y & \longrightarrow & Y & \xrightarrow{\{p^n\}_n} & \Pi_{n\geq 1} P^n Y \\ \| & & {\scriptstyle h}\uparrow{\scriptstyle \simeq} & & {\scriptstyle h}\uparrow{\scriptstyle \simeq} & & \| & & \\ \Omega Y & \xrightarrow{\overline{\Omega\{p^n\}_n}} & \overline{\Pi_{n\geq 1} \Omega P^n Y} & \longrightarrow & \overline{\Pi_{n\geq 1} \Omega P^n Y} & \longrightarrow & Y & & \end{array}$$

From universality (Proposition 5.14), we see that $\mathrm{Ph}^{\infty}_{\mathcal{F}}(X, Y) = *$ for any space $X$ if and only if the map $\overline{\Pi_{n\geq 1} \Omega P^n Y} \to Y$, which is itself a $\mathcal{F}$-phantom map with infinite Gray index, is trivial. By classical argument, this holds if and only if the map $\overline{\Omega\{p^n\}_n}$ admits a retraction. We then only have to show that if $\Omega Y$ is a retract of $\overline{\Pi_{n\geq 1} \Omega P^n Y}$, then $\overline{\Omega\{p^n\}_n}$ itself admits a retraction.

Suppose then there exists maps $\bar{i}$ and $\bar{r}$ such that $\bar{r} \circ \bar{i}$ is identity on $\Omega Y$:

$$\begin{array}{ccccc} \Omega Y & & \xrightarrow{\Omega p^n} & & \Omega P^n Y \\ {\scriptstyle \bar{i}}\downarrow & & & & \vdots\,{\scriptstyle i^n} \\ \overline{\Pi_{n\geq 1} \Omega P^n Y} & \xrightarrow[\simeq]{h} & \Pi_{n\geq 1} \Omega P^n Y & \xrightarrow[proj.]{} & \Omega P^n Y \\ {\scriptstyle \bar{r}}\downarrow & & & & \\ \Omega Y & & & & \end{array}$$

It is not to hard to see that $\Omega P^n Y \simeq P^{n-1}(\Omega Y)$, and that, through this equivalence, the map $\Omega p^n$ is identified as $p^{n-1}$. Using Point *3)* of Lemma 5.9, we deduce the existence of the dotted map $i^n$ making the diagram commutative. The map $h \circ \bar{i}$ is then the composition

$$\Omega Y \overset{\Omega\{p^n\}_n}{\longrightarrow} \Pi_{n\geq 1} \Omega P^n Y \overset{\Pi\, i^n}{\longrightarrow} \Pi_{n\geq 1} \Omega P^n Y,$$

and we deduce that the map $\bar{i}$ is the composition

$$\Omega Y \xrightarrow{\overline{\Omega\{p^n\}_n}} \overline{\Pi_{n\geq 1}\Omega P^n Y} \xrightarrow{\overline{\Pi\, i^n}} \overline{\Pi_{n\geq 1}\Omega P^n Y},$$

where $\overline{\Pi\, i^n}$ is a CW-approximation of $\Pi\, i^n$. Then $\bar{r}\circ\bar{i} \simeq \left(\bar{r}\circ\overline{\Pi\, i^n}\right)\circ\overline{\Omega\{p^n\}_n}$ is identity on $\Omega Y$, which means that $\overline{\Omega\{p^n\}_n}$ itself admits a retraction. □

*Proof of Proposition 2.10.* For any finite type rationally elliptic space $L$, the map $\pi_n(L)\to\pi_n(\widehat{L})$ induced in homotopy by the completion is an isomorphism for $n$ large enough, and then $p^n : Y \to P^n Y$ is an equivalence by Point *1)* of Lemma 5.9.

If $\Omega Y$ splits as stated, set $L_n$ as the product of all those $L_\alpha$ such that $L_\alpha \simeq P^n L_\alpha$, whereas $L_\alpha$ is not equivalent to $P^{n-1}L_\alpha$. If no such $L_\alpha$ exists, set $L_n = *$. We then clearly have a weak equivalence $\Omega Y \xrightarrow{\simeq} \Pi_{n\geq 0} L_n$. We deduce a weak equivalence $P^k\Omega Y \xrightarrow{\simeq} \Pi_{n\geq 0} P^k L_n$ for any $k$, and finally a weak equivalence

$$\overline{\Pi_{k\geq 1}\Omega P^k Y} \simeq \overline{\Pi_{k\geq 1} P^{k-1}\Omega Y} \xrightarrow{\simeq} \Pi_{k\geq 1}\Pi_{n\geq 0} P^{k-1}L_n = \Pi_{k\geq 0}\Pi_{n\geq 0} P^k L_n.$$

As $L_k \simeq P^k L_k$, then $L_k$ is a retract of $\Pi_{n\geq 0} P^k L_n$. We deduce that $\Pi_{k\geq 0} L_k$ is a retract of $\Pi_{k\geq 0}\Pi_{n\geq 0} P^k L_n$. Using the weak equivalences above, we see that $\Omega Y$ is a retract of $\overline{\Pi_{k\geq 1}\Omega P^k Y}$. □

## 6. The algebraic Gray filtration on SNT sets

We give below the proofs of main theorems of Section 3. Notice that the base point of $\mathrm{SNT}(X)$, described as $\varprojlim_n^1 \mathrm{Aut}(X^{(n)})$, is the homotopy type of the space $X$ itself. Therefore, when we say that some subset $A$ of $\mathrm{SNT}(X)$ is trivial, we mean in fact that $A = \{X\}$.

In what follows we suppose that $F$ is a functor, from the category of pointed homotopy types of connected CW-complexes to itself, and that there exists, for any space $X$ and any $n$, a natural equivalence $F(X^{(n)}) \simeq (FX)^{(n)}$. We then simply write $FX^{(n)}$.

In the following diagram between exact sequences, commutativity of the right square follows by naturality:

$$\begin{array}{ccccccc}
1 \longrightarrow & \mathrm{Aut}_F(X^{(n+1)}) & \longrightarrow & \mathrm{Aut}(X^{(n+1)}) & \longrightarrow & \mathrm{Aut}(FX^{(n+1)}) \\
& \vdots\downarrow & & \downarrow & & \downarrow \\
1 \longrightarrow & \mathrm{Aut}_F(X^{(n)}) & \longrightarrow & \mathrm{Aut}(X^{(n)}) & \longrightarrow & \mathrm{Aut}(FX^{(n)})
\end{array}$$

We then deduce the existence of the dotted arrow. This gives a map of tower $\{\mathrm{Aut}_F(X^{(n)})\}_n \to \{\mathrm{Aut}(X^{(n)})\}_n$, whose $\varprojlim_n^1$ we denote by

$$j_F : \varprojlim_n{}^1 \mathrm{Aut}_F(X^{(n)}) \to \varprojlim_n{}^1 \mathrm{Aut}(X^{(n)}) \cong \mathrm{SNT}(X).$$

*Proof of Theorem 3.2.* We owe the idea of the proof to a similar result of Félix-Thomas [3, Theorem 1], in a slightly different context. Instead of $F$, they consider the functor associating to a space $X$ the “homotopy group” $\bigoplus_{k\geq 1}\pi_k(X)$. Their proof however works as well here, and we content ourselves to describe some steps, without further details.

Let $P_F$ denotes the set of pairs $(Y, (f_n)_{n\geq 1})$, where $f_n : Y^{(n)} \to X^{(n)}$ are homotopy equivalences such that the following diagram commutes:

$$\begin{array}{ccc} FY^{(n+1)} & \xrightarrow{F(f_{n+1})} & FX^{(n+1)} \\ \downarrow & & \downarrow \\ FY^{(n)} & \xrightarrow{F(f_n)} & FX^{(n)} \end{array}$$

We introduce now an equivalence relation on $P_F$. Two pairs $(Y, (f_n)_{n\geq 1})$ and $(Z, (g_n)_{n\geq 1})$ in $P_F$ are said equivalent if there is an homotopy equivalence $\phi : Z \to Y$ such that, for any $n$, the maps $f_n \circ \phi \circ g_n^{-1}$ belongs to $\mathrm{Aut}_F(X^{(n)})$. There is then a well-defined map $P_F/\sim \xrightarrow{\rho} \mathrm{SNT}(X)$, associating to the class of $(Y, (f_n)_{n\geq 1})$ the homotopy type of $Y$.

The key step of the proof is then to show that the bijection $\theta : \mathrm{SNT}(X) \to \varprojlim_n^1 \mathrm{Aut}(X^{(n)})$ defined by Wilkerson [24] induces another map $\theta_F$ such that the following square commutes:

$$\begin{array}{ccc} P_X/\sim & \xrightarrow{\rho} & \mathrm{SNT}(X) \\ \downarrow{\scriptstyle\theta_F} & & \downarrow{\scriptstyle\theta} \\ \varprojlim_n^1 \mathrm{Aut}_F(X^{(n)}) & \xrightarrow[j_F]{} & \varprojlim_n^1 \mathrm{Aut}(X^{(n)}) \end{array}$$

Using Wilkerson techniques, we show also that $\theta_F$ is still a bijection.

As the image of $\rho$ is clearly $\mathrm{SNT}_F(X)$ (as defined in Section 3), the result follows. □

Notice that from the exact sequence of pointed sets

$$\varprojlim_n^1 \mathrm{Aut}_F(X^{(n)}) \xrightarrow{j_F} \varprojlim_n^1 \mathrm{Aut}(X^{(n)}) \to \varprojlim_n^1 \mathrm{Aut}(FX^{(n)}),$$

we deduce that $\mathrm{SNT}_F(X)$ belongs to the kernel of the natural map $\mathrm{SNT}(X) \to \mathrm{SNT}(FX)$, $Y \mapsto FY$. The following example shows that, in general, $\mathrm{SNT}_F(X)$ is a proper subset of that kernel.

**Example 6.1.** Let $X$ be a 1-connected finite type space, and let $F = -_{(0)}$ be the rationalisation functor. Then $\mathrm{Ker}\big(\mathrm{SNT}(X) \to \mathrm{SNT}(X_{(0)})\big)$ is $\mathrm{SNT}(X)$ itself, whereas $\mathrm{SNT}_{(0)}(X)$ is trivial.

*Proof.* The groups $\mathrm{Aut}_{(0)}(X^{(n)})$ are finite for all $n$, by [9, Corollary II.5.4]. Therefore, $\varprojlim_n^1 \mathrm{Aut}_{(0)}(X^{(n)}) \cong *$, and $\mathrm{SNT}_{(0)}(X)$ is trivial by Theorem 3.2. On another hand, $\mathrm{SNT}(X_{(0)}) = \{X_{(0)}\}$ by [24, Corollary II.b], and the example follows. □

It seems however difficult in general to characterize algebraically the "difference" between $\mathrm{SNT}_F(X)$ and $\mathrm{Ker}\big(\mathrm{SNT}(X) \to \mathrm{SNT}(FX)\big)$. One reason is that the image of $\mathrm{Aut}(X^{(n)}) \to \mathrm{Aut}(FX^{(n)})$ need not to be normal.

*Proof of Theorem 3.1.* Let $F$ be the functor $-^{(k)}$, the $k$-th Postnikov section. Applying $\varprojlim_n^1$ to the exact sequence $(n \geq k)$

$$1 \to \mathrm{Aut}_F(X^{(n)}) = \mathrm{Aut}_k(X^{(n)}) \to \mathrm{Aut}(X^{(n)}) \to G_k^n \to 1,$$

where $G_k^n$ stands for $\mathrm{Im}\big(\mathrm{Aut}(X^{(n)}) \to \mathrm{Aut}(X^{(k)})\big)$, we obtain the following exact sequence of pointed sets:

$$\varprojlim_n^1 \mathrm{Aut}_F(X^{(n)}) \xrightarrow{j_F} \varprojlim_n^1 \mathrm{Aut}(X^{(n)}) \overset{p^k}{\twoheadrightarrow} \varprojlim_n^1 G_k^n.$$

This identifies $\mathrm{SNT}^k(X) := \mathrm{Ker}\, p^k$ with the image of $j_F$. The result follows from Theorem 3.2 and the definition of $\mathrm{SNT}_F(X)$. □

*Proof of Theorem 3.3.* The finite type hypothesis on $X$ implies that the groups $\mathrm{Aut}(X^{(n)})$ are countable. Then Point *1)* follows directly from Theorem 2.15.

Let $G_n := \mathrm{Aut}(X^{(n)})$. Then the group $\mathrm{Aut}_k(X^{(k+1)})$ is $\mathrm{Ker}(G_{k+1} \to G_k)$, and the first statement of Point *2)* follows from Point *2)* of Theorem 4.4.

Look at now the fibration $X^{(k+1)} \xrightarrow{p^{(k)}} X^{(k)} \to K(\pi_{k+1}(X), k+2)$ where the right map is the $k$-th Postnikov invariant of $X$. Consider its associated exact sequence $[X^{(k+1)}, K(\pi_{k+1}(X), k+1)] \to [X^{(k+1)}, X^{(k+1)}] \xrightarrow{(p^{(k)})^*} [X^{(k+1)}, X^{(k)}]$. If $f \in \mathrm{Aut}_k(X^{(k+1)})$, then $(p^{(k)})^* \cdot f = p^{(k)} = (p^{(k)})^* \cdot Id$, which means that $f$ and $Id$ are in the same orbit through the action of $[X^{(k+1)}, K(\pi_{k+1}(X), k+1)] = H^{k+1}(X; \pi_{k+1}(X))$ on $[X^{(k+1)}, X^{(k+1)}]$. If this last cohomology group is finite, the same must be true for $\mathrm{Aut}_k(X^{(k+1)})$ (for a systematic study of the group $\mathrm{Aut}_k(X^{(k+1)})$, see [19]).

Point *3)* follows directly from Points *1)* and *2)*. □

The following example might seem a little bit artificial, but it has the advantage to illustrate drastically Remark 2.14. It is easily proved with Theorem 3.1.

**Example 6.2.** Let $X_1$ and $X_2$ be non equivalent spaces of the same $n$-type for all $n$. Then the set $L := \mathrm{SNT}(X_1) = \mathrm{SNT}(X_2)$ inherits in particular two different algebraic Gray filtrations: the first, whose subsets are denoted $L_1^k$, is inherited from the bijection $L \cong \varprojlim_n^1 \mathrm{Aut}(X_1^{(n)})$, whereas the second one, whose subsets are denoted $L_2^k$, is inherited from the bijection $L \cong \varprojlim_n^1 \mathrm{Aut}(X_2^{(n)})$.

This two filtrations are indeed in general distincts: suppose $X_2 \notin L_1^k$ for some $k$. This automatically holds for example if $X_1$ is rationally elliptic, by Point *3)* of Theorem 3.3. Then $X_1 \in L_1^\infty$ (as $X_1$ is the base point!), but $X_1 \notin L_2^k$.

We now decompose the proof of Theorem 3.4 in a sequence of three lemmas, each of which being a generalization, as well as a consequence, of similar lemmas of [14, 15]. In what follows, when we say that a group homomorphism $G \to H$ has a finite cokernel, it is meaned in fact that the image of $G$ has finite index in $H$.

**Lemma 6.3.** *If $X$ is as in Point 1) of Theorem 3.4, then, for all $n \geq k \geq 1$, the map $\mathrm{Aut}_k(X^{(n)}) \to \mathrm{Aut}_k(H^{\leq n}(X; \mathbb{Z}_\mathbb{P}))$ has finite kernel and finite cokernel. If $X$ is as in Point 2) of Theorem 3.4, then the map $\mathrm{Aut}_k(X^{(n)}) \to \mathrm{Aut}_k(\pi_{\leq n} X)$ has finite kernel and finite cokernel, for all $n \geq k \geq 1$.*

*Proof.* For the first claim, consider the following natural commutative diagram:

$$\begin{array}{ccccccc}
1 & \longrightarrow & \mathrm{Aut}_k(X^{(n)}) & \longrightarrow & \mathrm{Aut}(X^{(n)}) & \longrightarrow & \mathrm{Aut}(X^{(k)}) \\
 & & \downarrow & & \downarrow & & \downarrow \\
1 & \longrightarrow & \mathrm{Aut}_k(H^{\leq n}(X; \mathbb{Z}_\mathbb{P})) & \longrightarrow & \mathrm{Aut}(H^{\leq n}(X; \mathbb{Z}_\mathbb{P})) & \longrightarrow & \mathrm{Aut}(H^{\leq k}(X; \mathbb{Z}_\mathbb{P}))
\end{array}$$

where the rows are exact. By [14, Lemma 3.1.a], the two vertical maps on the right have finite kernel and cokernel. By chasing in this diagram, we then deduce the same is true for the left vertical map.

The second claim is proved similarly, using [15, Lemma 1]. □

The following key algebraic result is not stated this way in [14]. It is however not too hard to see that the following more general formulation stays valid.

**Lemma 6.4.** [14] *Let $\{G_n\}_n$ be a tower of countable groups. Suppose there exists $N$ such that each map $G_{k+1} \to G_k$, $k \geq N$, has a finite cokernel. Then $\varprojlim_n{}^1 G_n \cong *$ if and only if, for each $k \geq N$, the canonical map $\varprojlim_n G_n \to G_k$ has a finite cokernel.*

**Lemma 6.5.** *Let $X$ be 1-connected.*
*1) If the n-th Postnikov invariant $X^{(n)} \to K(\pi_{n+1}(X), n+2)$ is torsion, then $\mathrm{Aut}_k(X^{(n+1)}) \to \mathrm{Aut}_k(X^{(n)})$ has finite cokernel. If moreover $\pi_{n+1}(X) \otimes \mathbb{Q} \cong 0$, then it has also a finite kernel.*
*2) If $X$ is a co-$H_0$-space, then, for any $n$, $\mathrm{Aut}_k(X^{(n+1)}) \to \mathrm{Aut}_k(X^{(n)})$ has finite cokernel. Furthermore, if $H_{n+1}(X;\mathbb{Q}) \cong 0$ this map has also a finite kernel.*
*3) If $H^n(X;\mathbb{Q}) \cong 0$ for all $n > N$, then $\mathrm{Aut}_k(X^{(n+1)}) \to \mathrm{Aut}_k(X^{(n)})$ has a finite kernel and a finite cokernel for all $n \geq N$.*

*Proof.* We first show that for any space $X$, if, for some $n$, the map $\mathrm{Aut}(X^{(n+1)}) \to \mathrm{Aut}(X^{(n)})$ has finite kernel and/or cokernel, then the same must be true, for any $k$, for the map $\mathrm{Aut}_k(X^{(n+1)}) \to \mathrm{Aut}_k(X^{(n)})$. Indeed, if $n < k$, there is nothing to be proved. If $n \geq k$, this follows by chasing in the following natural commutative diagram, with exact rows:

$$\begin{array}{ccccccc}
1 & \longrightarrow & \mathrm{Aut}_k(X^{(n+1)}) & \longrightarrow & \mathrm{Aut}(X^{(n+1)}) & \longrightarrow & \mathrm{Aut}(X^{(k)}) \\
 & & \downarrow & & \downarrow & & \| \\
1 & \longrightarrow & \mathrm{Aut}_k(X^{(n)}) & \longrightarrow & \mathrm{Aut}(X^{(n)}) & \longrightarrow & \mathrm{Aut}(X^{(k)})
\end{array}$$

Then Points *1)* and *2)* follow from their corresponding statement in the case $k = 0$, to be found in [14, 15].

It remains to show Point *3)* in the case $k = 0$. As far as we know, this it not clearly stated in the literature, although our present proof uses techniques already present in [15].

We shall find a rational equivalence $W \xrightarrow{f} X$, where $W$ is a $N$-dimensional CW-complex (see for example the proof of Theorem 4 in [10]). Fix $n \geq N$. Let us define

$$\Delta(f^{(n)}) := \{(\alpha, \beta) \in \mathrm{Aut}(W^{(n)}) \times \mathrm{Aut}(X^{(n)}) \mid f^{(n)} \circ \alpha \simeq \beta \circ f^{(n)}\}.$$

Consider now the commutative diagram

$$\begin{array}{ccccc}
\mathrm{Aut}(W^{(n+1)}) & \longleftarrow & \Delta(f^{(n+1)}) & \longrightarrow & \mathrm{Aut}(X^{(n+1)}) \\
\cong \downarrow & & \downarrow & & \downarrow \\
\mathrm{Aut}(W^{(n)}) & \longleftarrow & \Delta(f^{(n)}) & \longrightarrow & \mathrm{Aut}(X^{(n)})
\end{array}$$

where the horizontal maps are given by the obvious projections. It follows from [25, Theorem 2.3] that these horizontal maps have finite kernel and cokernel. As $n$ is up the dimension of $W$, the left vertical map is in fact an isomorphism. By chasing in the diagram, we deduce the middle, and then the right vertical map has finite kernel and cokernel. □

*Proof of Theorem 3.4.* We only give the details for Point *1)*, the proof of Point *2)* being similar.

The equivalence of conditions *(ii)* and *(iii)* follows directly from Lemma 6.3.

The finite type hypothesis on $X$ implies that the groups $\mathrm{Aut}_k(X^{(n)})$ are countable, being subgroups of the countable groups $\mathrm{Aut}(X^{(n)})$. Now, the canonical surjection $\mathrm{Aut}(X) \twoheadrightarrow \varprojlim_n \mathrm{Aut}(X^{(n)})$ is easily seen to induce, for any $k$, a surjection $\mathrm{Aut}_k(X) \twoheadrightarrow \varprojlim_n \mathrm{Aut}_k(X^{(n)})$. Using Lemma 6.4, whose hypothesis is seen to be satisfied (with $N = 1$) thanks to Point *1)* of Lemma 6.5, we deduce that condition *(ii)* is equivalent to $\varprojlim_n^1 \mathrm{Aut}_k(X^{(n)}) \cong *$.

To conclude, we use now Point *1)* of Theorem 4.4. Indeed, with its notation, we have $K_k^n := \mathrm{Aut}_k(X^{(n)})$ and $L^k := \mathrm{SNT}^k(X)$. Therefore $\mathrm{SNT}^k(X)$ is trivial if and only if $\varprojlim_n^1 \mathrm{Aut}_k(X^{(n)}) \cong *$, and the result follows. □

*Proof of Example 3.5.* The space $BSU(m)$ is rationally equivalent to the product $K(\mathbb{Z}, 2) \times K(\mathbb{Z}, 4) \times \cdots \times K(\mathbb{Z}, 2m)$. Then $\pi_k(BSU(m)) \otimes \mathbb{Q} \cong 0$ for $k \geq 2m + 1$, and we deduce from Theorem 3.3 that $\mathrm{SNT}^{2m}(BSU(m))$ is trivial.

To show that $\mathrm{SNT}^{2m-1}(BSU(m))$ is non trivial, we use Point *1)* of Theorem 3.4, condition *(iii)*. The image of $\mathrm{Aut}(BSU(m)) \to \mathrm{Aut}(H^*(BSU(m); \mathbb{Z}))$ is finite (see the proof of [14, Theorem 4]). Then, for any $n \geq 2m$, the image of $\mathrm{Aut}_{2m-1}(BSU(m)) \to \mathrm{Aut}_{2m-1}(H^{\leq n}(BSU(m); \mathbb{Z}))$ is also clearly finite. It then suffices to show that $\mathrm{Aut}_{2m-1}(H^{\leq n}(BSU(m); \mathbb{Z}))$ is infinite. The cohomology algebra $H^*(BSU(m); \mathbb{Z})$ is the polynomial algebra $\mathbb{Z}[e_2, e_4, \ldots, e_{2m}]$ on generators $e_{2i}$ of degree $2i$. For $u \in \mathbb{Z}$, let $\phi_u$ be the automorphism sending $e_{2i}$ to itself for $1 \leq 1 \leq m - 1$, and sending $e_{2m}$ to $e_{2m} + u \cdot e_2 e_{2m-1}$. Now, if we associate to $u$ the restriction of $\phi_u$ to $H^{\leq n}(BSU(m); \mathbb{Z})$, we define an inclusion $\mathbb{Z} \hookrightarrow \mathrm{Aut}_{2m-1}(H^{\leq n}(BSU(m); \mathbb{Z}))$, and the result follows.

The assertions about $BSp(m)$ are showed similarly. □

To close this section, we provide the following weaker version of Theorem 3.4, available for spaces with a finiteness condition on their rational homotopy types.

**Theorem 6.6.** *Let $X$ be a 1-connected space, of finite type over some subring of the rationals. Suppose there exists $N$ such that $\pi_n(X) \otimes \mathbb{Q} \cong 0$ for any $n > N$, or such that $H_n(X; \mathbb{Q}) \cong 0$ for any $n > N$. Then $\mathrm{SNT}^k(X) = \{X\}$ if and only if one, and hence all, of the following equivalent conditions holds:*
*(i) The image of $\mathrm{Aut}_k(X) \to \mathrm{Aut}_k(X^{(N)})$ has finite index.*
*(ii) The image of $\mathrm{Aut}_k(X) \to \mathrm{Aut}_k(X^{(n)})$ has finite index, for some $n \geq N$.*
*(iii) The image of $\mathrm{Aut}_k(X) \to \mathrm{Aut}_k(X^{(n)})$ has finite index, for all $n \geq N$.*

*Proof.* As in the proof of Theorem 3.4, we see that $\mathrm{SNT}^k(X) = \{X\}$ if and only if $\varprojlim_n^1 \mathrm{Aut}_k(X^{(n)}) \cong *$.

By Point *1)* or Point *3)* of Lemma 6.5 , we see that all the maps $\mathrm{Aut}_k(X^{(n+1)}) \to \mathrm{Aut}_k(X^{(n)})$, $n \geq N$, have finite kernel and cokernel. This shows the equivalence of the three conditions. This shows also that the hypothesis of Lemma 6.4 are satisfied, and we then see that $\varprojlim_n^1 \mathrm{Aut}_k(X^{(n)}) \cong *$ if and only if condition *(iii)* is satisfied. □

## 7. From phantom maps to snt-theory

In that section we emphasize on the connection between phantom maps and SNT-theory, which is given by the maps

$$\mathrm{SNT}(Y \vee \Sigma X) \xleftarrow{\mathrm{Cof}} \mathrm{Ph}_{\mathcal{F}D}(X, Y) \xrightarrow{\mathrm{Fib}} \mathrm{SNT}(X \times \Omega Y)$$

associating to a $\mathcal{F}D$-phantom map its homotopy cofiber or its homotopy fiber. In what follows we suppose all our spaces to be 0-connected.

In [20] another connection is defined in an algebraic way. Let us define a group homomorphism $\psi_n : [X, \Omega Y^{(n)}] \to \mathrm{Aut}(X \times \Omega Y)^{(n-1)}$, associating to a map $f$ :

$X \to \Omega Y^{(n)} \simeq (\Omega Y)^{(n-1)}$ the homotopy equivalence

$$X^{(n-1)} \times (\Omega Y)^{(n-1)} \to X^{(n-1)} \times (\Omega Y)^{(n-1)}, \quad (x, \omega) \mapsto \big(x, \omega * f^{(n-1)}(x)\big)$$

where $*$ denotes here loop mutiplication.

The maps $\psi_n$ fit together in a map of towers $\{[X, \Omega Y^{(n)}]\}_n \to \{\mathrm{Aut}(X \times \Omega Y)^{(n-1)}\}_n$ and give rise to a map

$$\Psi : \mathrm{Ph}_{\mathcal{F}D}(X, Y) \cong \varprojlim_n{}^1 [X, \Omega Y^{(n)}] \longrightarrow \varprojlim_n{}^1 \mathrm{Aut}(X \times \Omega Y)^{(n-1)} \cong \mathrm{SNT}(X \times \Omega Y)$$

We don't detail the dual version, which starts from maps $\psi'_n : [\Sigma X, Y^{(n)}] \to \mathrm{Aut}(Y \vee \Sigma X)^{(n)}$ (there is no shift of degree here) and gives rise to a map

$$\Psi' : \mathrm{Ph}_{\mathcal{F}D}(X, Y) \cong \varprojlim_n{}^1 [\Sigma X, Y^{(n)}] \longrightarrow \varprojlim_n{}^1 \mathrm{Aut}(Y \vee \Sigma X)^{(n)} \cong \mathrm{SNT}(Y \vee \Sigma X)$$

It is said in [20] that it is not known if the maps $\Psi$ and $\Psi'$ are naturally equivalent, respectively, to Fib and Cof. The key result of that section is that indeed, they are.

**Theorem 7.1.** *The map $\Psi$ associates to a phantom map its homotopy fiber, whereas the map $\Psi'$ associates to a phantom map its homotopy cofiber.*

*Proof.* We only give the proof for the map $\Psi$. It may be dualized naturally to describe $\Psi'$.

Let $f : X \to Y$ be a $\mathcal{F}D$-phantom map, and let $F$ be the fiber of $f$. Working simplicially, we can construct the following strictly commutative diagram:

$$\begin{array}{ccccccccc}
F & \longrightarrow \cdots \longrightarrow & F_{n+1} & \xrightarrow{q_n} & F_n & \longrightarrow \cdots \xrightarrow{q_1} & F_1 \\
\downarrow & & \downarrow & & \downarrow & & \downarrow \\
X & = \cdots \longrightarrow & X & = & X & = \cdots = & X \\
{\scriptstyle f}\downarrow & & {\scriptstyle f^{n+1}}\downarrow & & {\scriptstyle f^{n}}\downarrow & & {\scriptstyle f^{1}}\downarrow \\
Y & \longrightarrow \cdots \longrightarrow & Y^{(n+1)} & \xrightarrow[\pi^n]{} & Y^{(n)} & \longrightarrow \cdots \longrightarrow & Y^{(1)}
\end{array}$$

where the bottom row is a Postnikov system for $Y$, and we define $F_n$ as the homotopy fiber of the map $f^n := X \to Y \to Y^{(n)}$. Then $F_n$ is the set of pairs $(x, \omega_n)$, where $x \in X$ and $\omega_n$ is a path in $Y^{(n)}$, starting at the base point, ending at $f^n(x)$. The natural maps $q_n$ are given by $(x, \omega_{n+1}) \mapsto \big(x, \pi^n \circ \omega_{n+1}\big)$.

We see that the maps in the upper row induce equivalences $F^{(n)} \simeq (F_{n+1})^{(n)}$. It follows they induce as well an equivalence from $F$ to the homotopy inverse limit of the tower $\{F_{n+1} \overset{q_n}{\to} F_n\}_n$.

Suppose now that, through the bijection $\mathrm{Ph}_{\mathcal{F}D}(X, Y) \cong \varprojlim_n{}^1 [X, \Omega Y^{(n)}]$, the map $f$ is described by some family $\{u_n : X \to \Omega Y^{(n)}\}_n$ . By construction of that bijection, that means we can choose homotopies $H^n : X \times I \to Y^{(n)}$, between $f^n$ and the trivial map ($H^n(x, 0) = f^n(x)$ and $H^n(x, 1) = *$), such that $u_n$ is the "difference" between $\pi^n \circ H^{n+1}$ and $H_n$. More precisely, we have the formula:

$$u_n : X \to \Omega Y^{(n)}, \quad x \mapsto \Big(\pi^n \circ \big(t \mapsto H^{n+1}(x, 1-t)\big)\Big) * \Big(t \mapsto H^n(x, t)\Big)$$

Now, homotopies $H^n$ can be used to construct homotopy equivalences:

$$v_n : F_n \xrightarrow{\simeq} X \times \Omega Y^{(n)}, \quad (x, \omega_n) \mapsto \big(x, \omega_n * (t \mapsto H^n(x, t))\big)$$

Let $h_n : X \times \Omega Y^{(n+1)} \to X \times \Omega Y^{(n)}$ by given by the formula $(x, \omega) \mapsto \big(x, \Omega\pi^n \cdot \omega * u_n(x)\big)$. Checking carefully, we see that we have construct a commutative diagram

$$\begin{array}{ccccc}
F_{n+1} & \xrightarrow[\simeq]{v_{n+1}} & X \times \Omega Y^{(n+1)} & \longrightarrow & (X \times \Omega Y)^{(n)} \\
& & & & \downarrow \\
{\scriptstyle q_n}\downarrow & & \downarrow{\scriptstyle h_n} & & (X \times \Omega Y)^{(n-1)} \\
& & & & {\scriptstyle \psi_n(u_n)}\downarrow \\
F_n & \xrightarrow[\simeq]{v_n} & X \times \Omega Y^{(n)} & \longrightarrow & (X \times \Omega Y)^{(n-1)}
\end{array}$$

where the unlabelled maps are induced by Postnikov sections, and $\psi_n$ are those maps involved in the definition of $\Psi$.

We deduce first that $F$ is homotopy equivalent to the homotopy inverse limit of the tower $\{X \times \Omega Y^{(n+1)} \overset{h_n}{\to} X \times \Omega Y^{(n)}\}_n$. Looking at homotopy groups, we then see that $F$ is also homotopy equivalent to the homotopy inverse limit of the tower

$$\big\{\, (X \times \Omega Y)^{(n)} \longrightarrow (X \times \Omega Y)^{(n-1)} \overset{\psi_n(u_n)}{\longrightarrow} (X \times \Omega Y)^{(n-1)} \,\big\}_n$$

By [24], that precisely means that $F$ is described as the family $\{\psi_n(u_n)\}_n$, through the bijection $\varprojlim_n^1 \operatorname{Aut}(X \times \Omega Y)^{(n-1)} \cong \operatorname{SNT}(X \times \Omega Y)$.

To summarize, we have shown that if $f \in \operatorname{Ph}_{\mathcal{F}D}(X, Y)$ is described by some family $\{u_n\}_n$, its homotopy fiber $F$ is described by the family $\{\psi_n(u_n)\}_n$. That exactly means that $\Psi(f) = F = \operatorname{Fib}(f)$. □

*Proof of Theorem 3.6.* By naturality of the construction of the algebraic Gray filtration, we see that a map $\varprojlim_n^1 G_n \to \varprojlim_n^1 G'_n$, which is the $\varprojlim_n^1$ of a map of towers $\{G_n \to G'_n\}_n$, must respect the filtrations. Then Theorem 3.6 is certainly true for $\Psi$ (Point *1*) and $\Psi'$ (Point *2*). Notice however the shift of degree for $\Psi$, due to the same shift for the maps $\psi_n : [X, \Omega Y^{(n)}] \to \operatorname{Aut}(X \times \Omega Y)^{(n-1)}$ involved in its definition.

Therefore Theorem 3.6 follows directly from Theorem 7.1. □

The following proposition generalizes Example 3.7, providing generic examples of spaces $Z$ for which the algebraic Gray filtration on $\operatorname{SNT}(Z)$ is highly non trivial.

**Proposition 7.2.** *Suppose that $X \simeq K(\mathbb{Z}, 2n)$, $n \geq 1$, and that $Y$, finite type with torsion fundamental group, is equivalent to $\Omega^j W$, for some $j \geq 0$ and some finite CW-complex $W$. Suppose furthermore that $H^n(X; \pi_{n+1}(Y) \otimes \mathbb{Q})$ is non zero for infinitely many $n$. Then there are infinitely many strict inclusions in the algebraic Gray filtration on* $\operatorname{SNT}(X \times \Omega Y)$.

*Proof.* Let $f : X \to Y$ be any phantom map. By [7, Lemma 3.3], the homotopy fiber $\operatorname{Fib}(f)$ of $f$ is not homotopy equivalent to $X \times \Omega Y$, provided $f$ is essential. Now, by Theorem 2.7, which is seen to apply thanks to Proposition 5.5, we deduce that there are, for infinitely many $n$, essential phantom maps $f_n : X \to Y$ of Gray index $n + 1$. As $\operatorname{Fib}(f_n) \in \operatorname{SNT}^n(X \times \Omega Y)$ (Theorem 3.6), we deduce that $\operatorname{SNT}^n(X \times \Omega Y)$ is non trivial for infinitely many $n$, and then that $\operatorname{SNT}^n(X \times \Omega Y)$ is non trivial for all $n < \infty$.

Suppose now only finitely many inclusions are strict in the Gray filtration on $\operatorname{SNT}(X \times \Omega Y)$. Then it would exist some integer $k$ such that $\operatorname{SNT}^k(X \times \Omega Y) = \operatorname{SNT}^\infty(X \times \Omega Y)$. By Point *1)* of Theorem 3.3, this would imply $\operatorname{SNT}^k(X \times \Omega Y)$ is trivial, a contradiction. □

*Proof of Example 3.8.* Let $W(X)/X$ be the cofiber of $\Phi_X : X \to W(X)$. As $\Phi_X \in \mathrm{Ph}^{\infty}_{\mathcal{F}D}(X, W(X))$, we deduce from Point *2)* of Theorem 3.6 that $W(X)/X$ belongs to $\mathrm{SNT}^{\infty}(W(X) \vee \Sigma X)$.

Suppose $X$ is finite type, and that its cohomology is not locally finite as a module over the Steenrod Algebra for some prime $p$. It is then proved in [18] that $\Phi_X$ is essential. Looking carefully the proof, we see it is in fact shown that $\Sigma X$ is not a retract of $W(X)/X$. In particular, $W(X)/X$ is not homotopy equivalent to $W(X) \vee \Sigma X$. □

Matematisk Institut, Universitetsparken 5, DK–2100 København

*E-mail address*: `ghienne@math.ku.dk`